\documentclass{article}

\usepackage{enumerate}
\usepackage{algorithm}
\usepackage{fullpage}
\usepackage{verbatim}
\usepackage{multicol}
\usepackage{graphicx}

 \usepackage{amsmath,amsfonts,amssymb, amsthm}
\usepackage{gensymb}
\usepackage{algorithm, algorithmicx, algpseudocode}

\def\omg{{D}}
\newcommand{\nRBF}{N_{RBF}}
\newcommand{\R}{\mathbb{R}}

\newcommand{\N}{\mathbb{N}}

\newcommand{\mbRd}{{\mbRn}}
\newcommand{\mbRn}{{\mathbb{R}^n}}

\def \zb{\mathbf{z}}
\def \xib{{\boldsymbol\xi}}
\newcommand{\fracOrder}{s}
\newcommand{\coloneqq}{:=}

\date{October 2022}

\title{Control of fractional diffusion problems\\ via dynamic programming equations}

\author{Alessandro Alla \thanks{Dipartimento di Scienze Molecolari e Nanosistemi, Università Ca' Foscari Venezia, Italy alessandro.alla@unive.it}, Marta D'Elia \thanks{Data Science and Computing, Sandia National Laboratories, Livermore, CA, USA mdelia@sandia.gov}, Christian Glusa \thanks{Center for Computing Research, Sandia National Laboratories, Albuquerque, NM, USA  caglussandia.gov}, Hugo Oliveira \thanks{Departamento de Matematica, PUC-Rio, Rio de Janeiro, Brazil, oliveira@mat.puc-rio.br}}
\begin{document}
\maketitle

\abstract{
  We explore the approximation of feedback control of integro-differential equations containing a fractional Laplacian term.
  To obtain feedback control for the state variable of this nonlocal equation we use the Hamilton--Jacobi--Bellman equation.
  It is well-known that this approach suffers from the curse of dimensionality, and to mitigate this problem we couple semi-Lagrangian schemes for the discretization of the dynamic programming principle with the use of Shepard approximation.
  This coupling enables approximation of high dimensional problems.
  Numerical convergence toward the solution of the continuous problem is provided together with linear and nonlinear examples.
  The robustness of the method with respect to disturbances of the system is illustrated by comparisons with an open-loop control approach.
}

\maketitle

\section{Introduction and motivation}\label{sec:introduction}
Nonlocal, integral models are powerful alternatives to classical partial differential equations (PDEs) to describe systems where small scale effects or interactions affect the global behavior, yet explicitly resolving local features is unfeasible due to computational budget or lack of local information.
Nonlocal models are characterized by integral operators that embed length scales in their definitions, allowing them to capture long-range space interactions.
Furthermore, the integral nature of such operators reduces the regularity requirements on the solutions.
Applications of interest span a large spectrum of scientific and engineering fields, including fracture mechanics \cite{Ha2011,Silling2000}, anomalous subsurface transport \cite{Benson2000,Gulian2021,Schumer2003,Schumer2001,suzuki2021fractional}, phase transitions \cite{Burkovska2021,Delgoshaie2015,Fife2003}, image processing \cite{Buades2010,DElia2021Imaging,Gilboa2007}, magnetohydrodynamics \cite{Schekochihin2008}, stochastic processes \cite{Burch2014,DElia2017,Meerschaert2012,MeKl00}, and turbulence \cite{DiLeoni2021,pang2020npinns,Pang2019fPINNs}.

Despite their improved accuracy, the usability of nonlocal equations is hindered by several modeling and computational challenges.
Modeling challenges include the lack of a complete nonlocal theory \cite{Defterli2015,DElia2020Helmholtz,DElia2020,d2021connections}, the treatment of nonlocal interfaces \cite{Alali2015,Capodaglio2020,Seleson2013,you2020asymptotically,yu2018partitioned} and the prescription of nonlocal volume constraints (the nonlocal counterpart of boundary conditions) \cite{DEliaNeumann2020,You_2019,trask2019asymptotically,yu2021asymptotically,Foss2021}.
Computational challenges are due to the integral nature of nonlocal operators that yields discretization matrices that feature a much larger bandwidth compared to the sparse matrices associated with PDEs \cite{DeliaDuEtAl2020_NumericalMethodsNonlocalFractionalModels,DEliaFEM2020}.

In the present work, we focus on fractional differential operators, as they form a well-studied sub-class of nonlocal models.
While these operators are almost as old as differential operators \cite{Gorenflo1997}, their usability has increased during the last decades thanks to progress in computational capabilities and thanks to a better understanding of their modeling power. The simplest form of a fractional operator is the fractional Laplacian; in its integral form, its action on a scalar function $y$ is defined as \cite{Gorenflo1997}
\begin{equation}\label{eq:fractional-laplacian}
  (-\Delta)^\fracOrder y(\xib) = C_{n,\fracOrder} \; {\rm p.v.}\int_\mbRn \frac{y(\xib)-y(\zb)}{\|\xib-\zb\|^{n+2\fracOrder}} d\zb,
\end{equation}
where $\fracOrder\in(0,1)$ is the {\it fractional order}, $n$ the spatial dimension, $C_{n,\fracOrder}$ a normalization constant and $\rm p.v.$ indicates the principal value.
The integral representation of the fractional Laplacian given in \eqref{eq:fractional-laplacian} is equivalent to the Fourier representation \cite{Valdinoci:2009}
\begin{equation}\label{E:definition_spectral}
(-\Delta)^{\fracOrder} y(\xib) =\mathcal{F}^{-1} ( \|\zb\|^{2\fracOrder} \mathcal{F}\{y\}(\zb) )(\xib),
\end{equation}
where $\mathcal{F}$ denotes the Fourier transform.
{If $\fracOrder=1$, the spectral operator \eqref{E:definition_spectral}
coincides with the usual PDE Laplacian $-\Delta$, whereas it reproduces the identity operator when $\fracOrder=0$.  In fact, it is {well known} that  $(-\Delta)^\fracOrder y(\xib) \rightarrow y(\xib)$ as $\fracOrder\to 0^+$ and $(-\Delta)^\fracOrder y(\xib) \rightarrow -\Delta y(\xib)$ as $\fracOrder\to  1^-$ for a regular function $y$; see e.g. \cite[Theorems~3 and~4]{Stinga:2019}}.

On a bounded Lipschitz domain $\omg\subset\mathbb R^d$, we define the {integral fractional Laplacian} to be the restriction of the full-space operator to functions satisfying a volume constraint on $\omg^{c}:=\mbRd\setminus\omg$. Here, for simplicity, we only consider the homogeneous case, i.e. $y=0$ in $\omg^{c}$. Then, the {fractional Poisson problem} on the bounded domain $D$ is given by
\begin{align}
  \left\{
  \begin{array}{rcll}
    (-\Delta)^{s}y(\xib) &=& b(\xib) &\text{for all }\xib\in\omg{,} \\
    y(\xib)&=&0 &\text{for all }\xib\in\omg^{c} ,
  \end{array}
  \right. \label{eq:fracPoisson}
\end{align}
where we have a given source term $b(\xib)$.
{Problem} \eqref{eq:fracPoisson} is a fractional analogue of the local Poisson problem for the classical Laplacian $\Delta$.

Accordingly, one can define a parabolic problem and add classical terms such as a nonlinear reaction.
The state equation that we consider in the present work is the time dependent, reaction-anomalous-diffusion equation given by
\begin{align}
  \left\{
  \begin{array}{rcll}
    \partial_t y(\xib,t) &=&-(-\Delta)^{\fracOrder}y(\xib,t) + F(y(\xib,t))+ b(\xib,t) &\text{for all }(\xib,t)\in\omg\times(0,\infty){,} \\
    y(\xib,t)&=&0 &\text{for all }(\xib,t)\in\omg^{c}\times(0,\infty), \\
    y(\xib,0) &=& x(\xib)&\text{for all } \xib\in\omg,
  \end{array}
  \right. \label{eq:fracReacDiff}
\end{align}
where $F(\cdot)$ is a reaction term, \(b\) a forcing and \(x\) the initial condition.

We are interested in the optimal control of equation \eqref{eq:fracReacDiff} with respect to a forcing term. While the PDE literature on control includes two main approaches, namely the open-loop and closed-loop approach, the nonlocal literature has only witnessed the development of the former. We refer the interested reader to \cite{d2014DistControl,DElia2016ParamControl,AO15,d2019priori,AntilKhatriEtAl2019_ExternalOptimalControlNonlocalPdes} for the the nonlocal Poisson equation and to \cite{burkovska2020,glusa2021error} for the fractional heat equation.

In the current work we, instead, aim at obtaining a control in closed form using the dynamic programming principle (DPP) and the Hamilton-Jacobi-Bellman (HJB) equations. This is, to the best of the authors' knowledge the first work in this direction. 
Feedback control has several advantages compared to open-loop control, since it explicitly depends on the state of the system, as opposed to open-loop where the
control depends only on time and the initial state. In the latter case, the fact that the control is not a function of the state at time $t$ has several disadvantages. For example, if the initial state changes, the optimal control has to be recomputed from scratch. Furthermore, open-loop control is not designed to handle modeling errors, exogenous disturbances, or other variations of system parameters.

 The DPP is an important tool in optimal control theory since it allows one to obtain feedback or closed-loop control. This is possible from the knowledge of the value function, i.e., the infimum of the cost functional of the control problem for every initial condition. Specifically, the value function is characterized as the unique viscosity solution of the HJB type PDE, see, e.g., \cite{BCD97}. 

Yet, closed-loop control also has some drawbacks. The main challenge of the DPP approach is the amount of information stored in the value function which makes the method suffer from the {\em curse of dimensionality}, i.e., the complexity greatly increases as the dimension does. In fact the control of a system of dimension $d$ requires the approximate solution of a $d$-dimensional PDE of HJB type. 
For high-dimensional dynamical systems, traditional methods such as the Value iteration or Policy iteration \cite{FF13} are not efficient to solve the HJB equation. In order to mitigate the curse of dimensionality, recently, several methods have been proposed. Here, we briefly recall some of them: model order reduction \cite{KVX04, AFV17}, tree structure algorithms \cite{AFS19,AS20}, spectral methods \cite{KK18}, max-plus algebra \cite{M07,M09}, neural networks \cite{DLM20, DM21}, tensor decomposition \cite{DKK21, OSS22} and sparse grids method \cite{BGGK13}.

In this work we use the approach recently proposed in \cite{AOS21}. The method relies on a high-dimensional scattered state-space grid that takes advantage of the dynamics of the problem and only populates certain regions using points from the discretized dynamics itself. The resulting unstructured grid is formed by nodes in the high-dimensional space. In order to directly apply a semi-Lagrangian scheme to approximate the value function, the Shepard method \cite{F07, F15} is used as the reconstruction tool. We also present a way to automate the selection of the shape parameter for the radial basis functions. The selection is driven by comparisons of the residual quantity. The coupling between HJB and Shepard approximation was already introduced by \cite{JS15} for a low-dimensional problem and equi-distributed grid. The algorithm proposed in \cite{AOS21} has also an automatic way to select the shape parameter in the chosen radial basis functions based on the residual of the HJB equation and error estimates for theoretical convergence of the method.

\paragraph{Paper outline}
In Section~\ref{sec:problemSetting} we describe the reaction diffusion state equation containing the fractional Laplacian as well as an associated optimal control problem.
We outline the discretization of the state equation using finite elements.
In Section~\ref{sec:HJBWithRadialBasisFunctions} we obtain the Hamilton-Jacobi-Bellman equation associated with the optimal control problem, and describe its discretization using radial basis functions.
Finally, in Section~\ref{sec:numerics}, we illustrate the potential of the approach via several numerical experiments.




\section{Problem settings}
\label{sec:problemSetting}
For the fractional Laplacian operator defined in \eqref{eq:fractional-laplacian} we let $D$ denote a bounded Lipschitz domain and we define the {integral fractional Laplacian} on a bounded domain to be the restriction of the full-space operator to functions satisfying a volume constraint on $\omg^{c}$. Here, for simplicity, we only consider the homogeneous case, i.e. $y=0$ in $\omg^{c}$.

We focus on the control of a general parabolic problem characterized by the following time-dependent, reaction-anomalous-diffusion equation for the variable $y:\mathbb{R}^n \times [0,\infty) \rightarrow \mathbb{R}$:
\begin{equation} \label{eq:time-dep-fractional-control}
	\left\{
	\begin{array}{rcll}
		\partial_t y(\xib,t) &=& -\alpha(-\Delta)^{\fracOrder} y(\xib,t) + F(y(\xib,t)) + & (\xib,t) \in \omg \times (0,\infty),\\
		\text{ } && \qquad \qquad \qquad \qquad  \quad  + b(\xib,t) + u(t)q(\xib) & \\
		y(\xib,t)&=&0 &(\xib,t)\in \omg^{c} \times (0,\infty) , \\
		y(\xib,0)&=&x(\xib) &\xib\in \omg.
	\end{array}
	\right.
\end{equation}
Here, $u:[0,\infty) \rightarrow U \subset \mathbb{R}$ is the control variable, where $U$ is a compact set, and $q: \omg \rightarrow \mathbb{R}$ is a known term; $F: \mathbb{R} \rightarrow \mathbb{R}$ is a nonlinear function, $b: D \times [0,\infty) \rightarrow \mathbb{R}$ is a forcing term and $\alpha > 0$ is a constant.

\noindent The general formulation of the cost functional that we intend to minimize is
\begin{align*}
  \mathcal{J}_{x}^{\infty}(y, u)
  &\coloneqq \int_0^\infty g(y(\eta),u(\eta)) e^{-\lambda \eta}d\eta,
\end{align*}
where, \(x\) is the initial condition, and $g: \mathbb{R}^d \times \mathbb{R}^m \rightarrow \mathbb{R}$ is the so-called running cost, assumed to be a bounded and Lipschitz continuous function in the first variable. The constant $\lambda > 0$ is a discount factor and the term $e^{-\lambda \eta}$ guarantees the convergence of the integral. The discount factor allows us to compare costs at different times by rescaling the costs at the initial time.
Then, the optimal control problem is formulated as
\begin{equation*}
	\min_{u \in U}	\mathcal{J}_{x}^{\infty}(y,u) \text{ s.t. } y(u) \text{ satisfying \eqref{eq:time-dep-fractional-control}}  ,
\end{equation*}
where we write $y = y(u)$ to emphasize the dependence of the solution on the control $u$.
 
We obtain the variational formulation of problem \eqref{eq:time-dep-fractional-control} by multiplying by a test function $w \in \mathcal{Y}$.
The function space $\mathcal{Y}$ is defined as
$$
\mathcal{Y} \coloneqq \{y \in H^{\fracOrder}(\mathbb{R}^n): y\|_{\omg^{c}}=0\},
$$
where, using the Fourier transform  $\mathcal{F}$, the fractional Sobolev space $H^s$ is defined as \cite{mclean2000strongly}
\begin{align*}	H^{\fracOrder}\left(\mathbb{R}^n\right)&:=\left\{y\in L^{2}\left(\mathbb{R}^n  \right) : \int_{\mathbb{R}^n}(1+\|\xib\|^{2 \fracOrder})\|\mathcal{F}y(\xib) \|^2 d\xib < \infty \right\}.
\end{align*}
We equip $\mathcal Y$ with the following norm
\begin{align*}
\|y\|_{\mathcal{Y}}^{2}&=\|y\|_{H^{\fracOrder}(\mathbb{R}^n)}^{2}= \|y\|_{L^{2}\left(\omg\right)}^{2} + \int_{\mathbb{R}^n} \int_{\mathbb{R}^n} \frac{\left(y(\xib)-y(\zb)\right)^{2}}{\|\xib-\zb\|^{n+2\fracOrder}}d\xib d\zb.
\end{align*}
Hence, the variational formulation of the parabolic problem \eqref{eq:time-dep-fractional-control} is given by
\begin{align}
	\label{variational}
	\text{Find } y  \in \mathbb{V} &\text{ such that } \forall \;w\in \mathcal{Y} \nonumber\\
	\langle \partial_t y,w \rangle = -\alpha a(y,w) + & \langle F(y),w\rangle + \langle b,w\rangle + u\langle q,w\rangle,
\end{align}
where the space $\mathbb{V}$ is defined as
$$
\mathbb{V} \coloneqq \{y \in L^{2}((0,\infty);\mathcal{Y}): \partial_t y \in L^{2}((0,\infty);\mathcal{Y'})\},
$$
$\mathcal{Y'}$ is the dual space of $\mathcal{Y}$, and \(\langle \cdot,\cdot \rangle\) is the duality pairing between \(\mathcal{Y}\) and its dual.
The term $a(y,w)$ is obtained by applying integration by parts \cite{DElia2020} to the fractional Laplacian \eqref{eq:fractional-laplacian}, resulting in
\begin{align*}
	a(y,w)
	&\coloneqq
	 \frac{C_{n,s}}{2}\int_{\mathbb{R}^n} \int_{\mathbb{R}^n} (y(\xib)-y(\zb)) (w(\xib)-w(\zb)) \|\xib-\zb\|^{-n-2\fracOrder} d\zb d\xib.
\end{align*}

We next describe the finite element discretization of the variational formulation \eqref{variational}.
For details on issues such as quadrature rules suitable for resolving the singularity of the kernel function, see \cite{AinsworthGlusa2018_TowardsEfficientFiniteElement,AinsworthGlusa2017_AspectsAdaptiveFiniteElement}.

Let $\mathcal{G} = \{ \mathcal{K} \}$ be a conforming partition of $\overline \omg$ into simplices $\mathcal{K}$ with size $h_\mathcal{K} = \operatorname{diam}(\mathcal{K})$, and set $h_{\mathcal{G}} = \max \limits_{\mathcal{K} \in \mathcal{G}} h_\mathcal{K}$.
Given $\mathcal{G}$, we define the finite element space of continuous piecewise polynomials of degree one as
\begin{equation}
	\mathcal{Y}(\mathcal{G}) = \left\{ w_{\mathcal{G}} \in C^0( \overline \omg): {w_{\mathcal{G}}}_{\|\mathcal{K}} \in \mathbb{P}_1(\mathcal{K}),  \forall \mathcal{K} \in \mathcal{G}, \ w_{\mathcal{G}} = 0 \textrm{ on } \partial D \right\},
\end{equation}
where $\mathbb{P}_1 (\mathcal{K})$ is the space of linear functions on $\mathcal{K}$. We point out that $\mathcal{Y}(\mathcal{G}) \subset \mathcal{Y}$.
Let $\{\phi_i\}_{i=1}^{d}\subset \mathcal{Y}(\mathcal{G})$ be the nodal hat functions that span the space. We search for a function $y_\mathcal{G} \in \mathcal{Y}(\mathcal{G})$ which can be described by the ansatz $y_\mathcal{G}(\xib,t)=\sum_{i=1}^{d} y_{\mathcal{G},i}(t) \phi_i (\xib)$, i.e. for every fixed time $t$ the function $y_\mathcal{G}$ is a continuous piecewise linear function with time dependent nodal values $y_i (t)$. Thus, using this ansatz in \eqref{variational}, we obtain the terms
\begin{align*}
	a(y_\mathcal{G},\phi_j) &=\sum_{i=1}^{d} y_{\mathcal{G},i} a(\phi_i,\phi_j), &
	(\dot{y_\mathcal{G}},\phi_j) &=\sum_{i=1}^{d} \dot{y_{\mathcal{G},i}} (\phi_i,\phi_j),
\end{align*}
where $\dot{y} = \frac{dy}{dt}$. Using these spatially discretized entities and treating the other terms accordingly (note that the initial condition can be taken to be $x_{\mathcal{G}}(\xib) = \sum_{i=1}^{d} x_{\mathcal{G},i} \phi_i (\xib)$ after projection or interpolation), we have the semi-discrete version of equation \eqref{eq:time-dep-fractional-control} given by
\begin{equation} \label{fin_element}
	\begin{cases}
		M\dot{y_{\mathcal{G}}}(t)= -\alpha Ay_{\mathcal{G}}(t) + \mathbf{F}(y_{\mathcal{G}}(t)) + B(t) + u(t)Q, \quad & t\in(0,\infty),\\
		y_{\mathcal{G}}(0)=x_{\mathcal{G}}.
	\end{cases}
      \end{equation}
We note that by abuse of notation we have identified finite element functions such as \(y_\mathcal{G}=\sum_{i=1}^{d} y_{\mathcal{G},i} \phi_i\) with its vector of coefficients \(\{y_{\mathcal{G},i}\}_{i=1}^{d}\).
Here, the matrices $A,M \in \mathbb{R}^{d \times d}$ are constructed with entries
$	A_{ij} = a(\phi_i,\phi_j), 
	M_{ij} = (\phi_i,\phi_j).$
The other elements are $d$-dimensional vectors with entries given by $\mathbf{F}_i (y_{\mathcal{G}}) = \langle F(y_{\mathcal{G}}), \phi_{i} \rangle$, $B_i (t) = \langle b(t,\cdot),\phi_i\rangle$ and $Q_i  = \langle q,\phi_i\rangle$.
Since in what follows, we only operate at the discrete level, we drop the subscript \(\mathcal{G}\) and we rewrite the semi-discrete equation as
\begin{equation} 
	\begin{cases}
		\dot{y} (t)= M^{-1}(-\alpha Ay(t) + \mathbf{F}(y(t)) + B(t) + u(t)Q), \quad & t\in(0,\infty),\\\label{dyn_frac}
		y(0)=x.
	\end{cases}
\end{equation}
For the sake of clarity and consistency with the following section (see equation \eqref{dyn}), we denote the right-hand side of the dynamical system in equation \eqref{dyn_frac} as follows
\begin{equation}\label{disc_non}
f(y(t),u(t)) = M^{-1}(-\alpha Ay(t) + \mathbf{F}(y(t)) + B(t) + u(t)Q).
\end{equation}

\section{The HJB approach with radial basis functions}
\label{sec:HJBWithRadialBasisFunctions}

We recall the main results on infinite-horizon control problems solved by means of dynamic programming equations. For a complete description we refer to, e.g., the manuscripts \cite{BCD97,FF13}. We consider the dynamical system
\begin{equation}  
	\begin{cases}
		\dot{y} (t)=f(y(t),u(t)), \text{ } t \in (0,\infty),\\ \label{dyn}
		y(0)=x \text{}\in \mathbb{R}^d,
	\end{cases}
\end{equation}
where $y: [0,\infty) \rightarrow \mathbb{R}^d$ is the state variable, and $u\in \mathcal{U} \coloneqq \{ u: [0,\infty) \rightarrow U, \text{ measurable}\}$ is the control variable. Here, $\mathcal{U}$ is the set of admissible controls, $U\subset\R^m$ is a compact set, and $f:\mathbb{R}^d \times \mathbb{R}^m \rightarrow \mathbb{R}^d $ is a Lipschitz continuous function with respect to the first variable with constant $L_f>0$. 
Under such hypothesis, \eqref{dyn} has a unique solution \cite{FF13}. 

As anticipated in the previous section, the cost functional $\mathcal{J}: \mathcal{U} \rightarrow \mathbb{R}$ is given by
\begin{equation}
 \label{costfunctional}
		\mathcal{J}_x(y, u) = \int_0^\infty g(y(\eta),u(\eta)) e^{-\lambda \eta}d\eta,
\end{equation}
with running cost \(g\) and discount factor \(\lambda\). Thus, the optimal control problem reads
\begin{equation}
	\min_{u \in \mathcal{U}} \mathcal{J}_x(y, u),
\end{equation}
with $y$ being a trajectory that solves \eqref{dyn} corresponding to the initial point $x$ and control $u$. The use of the index $x$ emphasizes the dependence on the initial condition. Since we aim at obtaining the control in a feedback form, i.e., as a function of the current position of the state, we define the value function as
\begin{equation}
	\label{vf}
	v(x) \coloneqq  \inf_{u \in \mathcal{U}} \mathcal{J}_x(y, u).
\end{equation}
It is possible to characterize the value function in terms of the DPP as follows
\begin{equation}
	\label{DPP}
	v(x) = \inf_{u \in \mathcal{U}} \bigg\{\int_0^\tau g(y(\eta),u(\eta)) e^{-\lambda \eta}d\eta + e^{-\lambda \tau}\ v(y_{x}(\tau))\bigg\} \quad  \forall x \in \mathbb{R}^d,\tau>0.
\end{equation}
Then, from \eqref{DPP}, the Hamilton-Jacobi-Bellman (HJB) equations corresponding to the infinite horizon problem are given by
\begin{equation}
	\label{HJB}
	\lambda v(x) + \sup_{u \in U} \{ -f(x,u) \cdot \nabla v(x) - g(x,u)\}=0, \qquad x\in\mathbb{R}^d,
\end{equation}
where $\nabla v$ is the gradient of $v$. The HJB is a further characterization of the value function by means of a degenerate PDE whose solution has to be thought of in the viscous sense \cite{BCD97}.
Thus, by solving \eqref{HJB}, we obtain the optimal feedback control $u^{*}(x)$ as follows

\begin{equation}
	\label{FDB}
	u^{*}(x) = \operatorname*{ arg\,max}_{u \in U} \{ -f(x,u) \cdot \nabla v(x) - g(x,u)\},\qquad x\in\mathbb{R}^d.
\end{equation}

To approximate the solution of \eqref{HJB}, the semi-Lagrange discretization of \eqref{DPP} by means of Shepard approximation was first introduced in \cite{JS15} and then extended to higher-dimensional problems in \cite{AOS21}. Such discretization of \eqref{DPP} reads
\begin{equation}
	\label{fully-discrete}
	V_j = [W_\sigma(V)]_{j}\coloneqq \min_{u \in U} \bigg\{  \Delta t\,g(x_j,u) + (1-\Delta t \lambda) S^\sigma[V](x_j + \Delta t\, f(x_j,u))\bigg\},
\end{equation}
where $V_j\approx v(x_j)$, $\Delta t>0$ is the chosen temporal step size, and $x_j$ are nodes of the grid $X=\{x_1,\ldots, x_{N_{grid}}\}\subset\Omega\subset\R^d$ which we describe later on. The set $U$ is also discretized by setting the control $u(t)=u_k \in U$ for $t \in [t_k, t_{k+1})$, by approximating the control by a piece-wise constant function in time. The minimum in \eqref{fully-discrete} is computed by comparison over the discrete set of controls. 

The full discretization 
of equation \eqref{DPP} is obtained starting from a classical approach where we replace the interpolation routine on a structured grid (see e.g. \cite{FF13}) with the Shepard approximation operator $S^\sigma$ on an unstructured grid.
We represent the Shepard approximant as the following operator
\begin{equation}\label{eq:shepard_basis}
 S^\sigma[f](x) \coloneqq \sum_{i=1}^{\nRBF} f(x_i)\psi_i^\sigma (x), \quad
\psi^\sigma_i(x)\coloneqq  \frac{\varphi^\sigma(\|x-x_i\|_2)}{\sum_{j=1}^{\nRBF} \varphi^\sigma(\|x-x_j\|_2)}, \;\; 1\leq i\leq \nRBF.
\end{equation}
where $\sigma>0$ is the shape parameter, and $\nRBF$ the number of RBFs $\varphi$.

Here, the function $\varphi: \mathbb R^d \rightarrow \mathbb R$ is a positive definite radially invariant function $\varphi (x) = \varphi (\|x\|_2)$.
We refer to e.g. \cite[Chapter 23]{F07} for a complete review of Shepard approximation scheme. The shape parameter $\sigma$ is used to tune the radial function such that $\varphi^{\sigma} (x) = \varphi (\sigma \|x\|_2)$. The effect of the shape parameter in the quality of Shepard approximation is discussed in e.g. \cite{F07, AOS21}.

We point out that each $\psi_i(x)$ is compactly supported in $B(x_i, 1/\sigma)\subset\Omega$ and nonnegative, and that the weights form a partition of unity, 
i.e., $\sum_{i=1}^{\nRBF} \psi_i^\sigma(x) = 1$ for all $x\in\Omega_{X, \sigma}$ with
\begin{equation}\label{eq:omega_z_sigma}
\Omega_{X, \sigma}\coloneqq \bigcup_{x\in X} B(x, 1/\sigma)\subset\R^d.
\end{equation} 
This implies that $S^\sigma[f](x)$ is a 
convex combination of the function values. Moreover, the compact support of the weights leads to a sparse distance vector $D\in\R^{\nRBF}$ since we compute only the entries $D_i$ such that $\|x-x_i\|\leq 1/\sigma$. 

The full approximation scheme of the value function is known as the Value Iteration (VI) method, and it is obtained by a fixed-point iteration of \eqref{fully-discrete}, i.e.,
\begin{equation}
\label{fixedpoint}
    V^{k+1} = W_\sigma(V^{k}),\quad k =0,1,\ldots.
\end{equation}
To ensure convergence of the scheme it is necessary that $W_\sigma$ is a contraction which is guaranteed to be the case since $S^\sigma$ in
\eqref{eq:shepard_basis} has unit norm if $\Delta t \in \left(0,1/\lambda\right]$ (see \cite[Lemma 2]{JS15}). The approximation of the value function is then used to compute an approximation of the feedback control as
\begin{equation}
	u^{*}(x) = \operatorname*{ arg\,min}_{u \in U} \{\Delta t g(x,u) + (1-\lambda \Delta t )S^\sigma[V](x+\Delta t f(x,u)) \},
\end{equation}
with $x$ being the current state of the discrete system. This also allows us to reconstruct the optimal trajectory $y^{*}$.

\paragraph{\bf Mesh generation}
The key point suggested in \cite{AOS21} is to localize the domain of the problem. Thus, we follow the evolution of the system \eqref{dyn} which provides itself an indication of the regions of interest within the
domain.

We point out that for the purpose of numerical computations we consider only finite horizons $t\in[0, T]$ with a given $T>0$ large enough to simulate the infinite-horizon problem.
To define the dynamics-dependent grid we fix a time step $\overline{\Delta t}>0$, a maximum number $\overline K\in\N$ of discrete times and, for $\overline 
L, \overline M >0$, some initial conditions of interest and a discretization of the control space, i.e., 
\begin{equation*}
\overline{X}\coloneqq \{\bar x_1, \bar x_2, \ldots, \bar x_{\overline L}\} \subset \Omega, \quad
\overline{U}\coloneqq \{\bar u_1,\bar u_2, \ldots, \bar u_{\overline{M}}\}\subset U.
\end{equation*}
Observe that all these parameters do not need to coincide with the ones used in the solution of the VI \eqref{fully-discrete}, but they are only used to construct the grid. Moreover, in general we use $\overline{\Delta t} > \Delta t$ and $\overline{M}< M$, i.e., the discretization 
used to construct the mesh is coarser than the one used to solve the control problem.

For a given pair of initial condition $\bar x_i\in \overline{X}$ and control $\bar u_j\in \overline{U}$ we solve numerically \eqref{dyn} to obtain trajectories
\begin{align}\label{grid-dyn}
x_{i,j}^{k+1} & = x_{i,j}^k + \overline{\Delta t}\, f(x_{i,j}^k, \bar u_j), \quad k=1,\ldots,\bar K-1,\\
x_{i,j}^1& = \bar x_i,\notag
\end{align}
such that $x_{i,j}^k$
is an approximation of the state variable with initial condition $\bar x_i$, constant control $u(t)=\bar u_j$, at time $k
\overline{\Delta t}$. For each pair $\left(\bar x_i, \bar u_j\right)$ we obtain the set $X\left(\bar x_i, \bar u_j\right) 
\coloneqq \{x_{i,j}^1,\ldots,x_{i,j}^{\bar K}\}$ containing the discrete trajectory, so that the mesh is defined as
\begin{equation}\label{eq:mesh}
X\coloneqq X(\overline X, \overline U, \overline \Delta t, \overline K)\coloneqq \bigcup\limits_{i=1}^{\bar L} \bigcup\limits_{j=1}^{\bar M} X\left(\bar x_i, \bar u_j\right).
\end{equation}
 This choice of the grid is particularly well suited for the problem under consideration, as it does not aim at filling the space $\Omega$, but instead it 
provides points along trajectories of interest. In this view, the values of   $\overline X, \overline U, \overline \Delta t, \overline K$ should be 
chosen such that $X$ contains points that are suitably close to the points of interest for the solution of the control problem. 
To measure the ``dispersion'' of points in $X$ we make use of the \emph{separation distance} $h$
\begin{equation*}
h = h_{X} \coloneqq\min_{x_i \neq x_j\in X}\|x_i - x_j \|.
\end{equation*} 
The separation distance
quantifies the minimal separation between different approximation points; as such, we use $h$ as the measure of dispersion in the numerical tests. A more detailed discussion of these quantities can be found in e.g. \cite{AOS21,F07,JS15}.
The quality of the Shepard approximation strongly depends on the choice of the shape parameter $\sigma>0$ in the operator \eqref{eq:shepard_basis}. To get a suitable scale for the values of $\sigma$, we parametrize it in terms of the grid $X$: we set $\sigma\coloneqq  \theta/h_{\Omega, X}$ for a given $\theta >0$.
Here, $\theta \in \mathcal{P}$ where $\mathcal{P}$ is a positive closed interval. In this work we follow the method proposed in \cite{AOS21}, where the selection of the ideal parameter $\sigma$ is based on the minimization of a problem-specific indicator, namely the residual 
$R(\sigma)\coloneqq R(V_{\sigma})$. Here, $V_{\sigma}$ is the discrete value function related to the choice $\sigma$ in the Shepard approximation.
Assuming that the value iteration with parameter $\sigma$ converges to
$V_{\sigma}$, we define the residual as 
\begin{equation}\label{residual}
	R(\sigma)\coloneqq  \|V_{\sigma} - W_\sigma(V_{\sigma}) \|_\infty,
\end{equation}
and we choose the shape parameter that minimizes this quantity with respect to $\sigma$. Since the minimization of the residual depends on $\sigma = \theta/ h_X$, we set an admissible set of parameters $\mathcal{P} \coloneqq [\theta_{\min}, \theta_{\max}]\subset \mathbb{R}^{+}$, the parameter is thus chosen by solving 
the optimization problem	
\begin{equation}\label{opt-shape}
	\bar\theta\coloneqq  \operatorname*{ arg\,min}\limits_{\theta\in\mathcal{P}} R(\theta/h_X) =  \operatorname*{ arg\,min}\limits_{\theta\in\mathcal{P}}\|V_{\theta/h_X} - W_{\theta/h_X}(V_{\theta/h_X}) 
	\|_\infty,
\end{equation} 
by discretizing the set $\mathcal{P}$ as $\{\theta_1, \cdots \theta_{N_p}\}\subset \mathcal{P}$ and computing all value functions for all $\theta_i$, $i=1,\cdots,N_p$.

\paragraph{\bf Error estimates} 
We adapt the classical convergence theory that is used to prove rates of convergence for the VI when linear interpolation is used. In 
particular, the following argument follows the discussion in \cite[Section 8.4.1]{FF13}. Provided that
\begin{itemize}
\item $v$ is Lipschitz continuous, with constant $L_v$ for each compact subset $\tilde \Omega\subset \Omega$,
\item $f$ is uniformly bounded, $\sup\limits_{x\in\Omega, u\in U} \|f(x, u)\| \leq M_f$ with $M_f>0$,
\item $U$ is convex, $f(x, u)$ is linear in $u$,
\item $g(\cdot, u)$ is Lipschitz continuous and $g(x, \cdot)$ is convex,
\item $\|f(x, u)-f(x', u)\| \leq L_{x} \|x-x'\|\;\; \forall x, x'\in \Omega,\;\; u\in U,$ $L_x>0$, and
\item $\|f(x, u)-f(x, u')\| \leq L_{u} \|u-u'\|\;\; \forall x\in \Omega,\;\; u, u'\in U$, $L_u>0$,
\end{itemize}
then, there exist $C,C'>0$ such that the following estimate for the value function holds true:

\begin{align*}
    \|v - V\|_{\infty,\tilde\Omega}\leq C'\Delta t +  \frac{C L_v}{\lambda} \frac{e^{\bar K \overline \Delta t L_x}}{\Delta t} \left(\sup\limits_{\Delta t\in \tilde T}\|\Delta t - \overline{\Delta t}\| \bar K
M_f \right.\\
\qquad\qquad\left.+ \sup\limits_{x_0\in \tilde X}\min\limits_{\bar x\in \bar X} \|\bar x- x_0\| + \bar K \overline{\Delta t} L_u \sup\limits_{u\in \tilde U}\min_{\bar u\in \bar U} 
\|\bar u- u\|\right) 
\end{align*}
where 
\begin{equation}\label{eq:set_omega_tilda}
\tilde \Omega\coloneqq  \tilde \Omega(\tilde X, \tilde U, \tilde T)\coloneqq  \left\{x\coloneqq  x^k(x_0, u, \Delta t): x_0\in \tilde X, u\in \tilde U, \Delta 
t\in \tilde T, k\leq \bar K\right\}.
\end{equation}
Here, $x\coloneqq x^k(x_0, u, \Delta t)\in\Omega$ is a point on a discrete trajectory with initial point $x_0\in\Omega$, control $u\in U$, time step $\Delta t>0$, 
and time instant $k\in \N$, $k\leq \bar K$. We point out that the $f$ defined in \eqref{disc_non} is uniformly bounded and Lipsichitz continuous with respect to both variables as required.

\section{Numerical Tests}
\label{sec:numerics}
In order to illustrate the DP approach to control fractional diffusion problems we present three tests. The first is a parabolic problem with known analytical solution. In this case, we show a comparison between the analytical solution and the one obtained by the DP approach. The same equation is also used to show the effectiveness of the feedback control when dealing with disturbances of the system. We also discuss the numerical convergence of our approach towards the control of the continuous fractional PDE.

The second and third tests are performed using exclusively the DP approach. The second test is a modification of previous problem where we control the equation on a target subset of the domain $D \subset \mathbb{R}$. The third test is a nonlinear state equation with a quadratic cost functional.

In all our tests we use Wendland RBFs \cite{W95, F07} constructed according to the dimension of each problem obtained using the general formula proposed in \cite{F98} $$\varphi^{\sigma}(r)= \max \{0, (1-\sigma r)^{\ell + 2} ((\ell^2 + 4\ell + 3)\sigma^2 r^2  + (3\ell + 6) \sigma r + 3)\}, \;\; r\coloneqq \|x\|_2,$$
where $\ell = \lfloor \frac{d}{2} \rfloor + 3$. The numerical simulations reported in this paper are performed on a laptop with one CPU Intel Core i7-2.2 GHz and 16GB RAM. The codes are written in MATLAB.   

\subsection{Test 1: Linear example with exact solution}\label{sec:example}

In this test, with the goal of comparing the solution obtained with the DP approach with a known solution, we consider the following manufactured problem.

Let \(\gamma>0\), \(\lambda>0\), \(T_{0}>0\) and set
\begin{align}
	\label{terms}
	\omg &= (-1,1), \nonumber\\
	\tilde{b}(\xib) &= 2^{2\fracOrder}\Gamma\left(1+\fracOrder\right) \frac{\Gamma(\fracOrder+1/2)}{\Gamma(1/2)} \sqrt{\frac{\Gamma(2\fracOrder+3/2)}{\Gamma(2\fracOrder+1)\Gamma(1/2)}}, \nonumber \\
	q(\xib) &= \sqrt{\frac{\Gamma(2\fracOrder+3/2)}{\Gamma(2\fracOrder+1)\Gamma(1/2)}} \left(1-\xib^{2}\right)^{\fracOrder}_{+},\\
	\phi(t)&= \cos(t), \nonumber \\
	\kappa(t)&= (T_{0}-t)^{2}\chi_{\{t\leq T_{0}\}}, \nonumber \\
	u_{d}(t) &= \operatorname{proj}_{U}(\kappa(t)), \nonumber \\
	y_{d}(\xib,t) &=\phi(t)q(\xib) - \gamma \kappa'(t)q(\xib) + \gamma \kappa(t) \tilde{b}(\xib) + \lambda\gamma \kappa(t)q(\xib), \nonumber \\
	b(\xib,t) &= \phi'(t)q(\xib) + \phi(t)\tilde{b}(\xib) - u_{d}(t)q(\xib). \nonumber
\end{align}

We minimize
\begin{align}
	\label{cost_funct1}	
	\mathcal{J}_{x}(y, u) &= \frac{1}{2} \int_0^{\infty} (\|y(\cdot,\eta)-y_{d}(\cdot,\eta)\|_{L^{2}(\omg)}^{2}  + \gamma \|q(\cdot)u(\eta)\|_{L^{2}(\omg)}^2) e^{-\lambda \eta} d\eta
\end{align}
subject to
\begin{equation} \label{eq:test1}
	\left\{
	\begin{array}{rll}
		\partial_t y(\xib,t) + (-\Delta)^{\fracOrder}y(\xib,t) &= b(\xib,t) + u(t)q(\xib) &(\xib,t) \in \omg \times (0,\infty),  \\
		y(\xib,t)&=0 &(\xib,t)\in \omg^{c} \times (0,\infty) , \\
		y(\xib,0)&=x(\xib) &\xib\in \omg.
	\end{array}
	\right.
\end{equation}

As shown in \ref{sec:deriv-example} the solution for \(x=q\) is given by
\begin{align}
	\label{eq:test1_optimal}
	y^{*}(\xib,t)=\phi(t)q(\xib), \qquad
	u^{*}(t)=u_{d}(t).
\end{align}

\paragraph{Change of variables}

To fit the cost functional into the setting of the DPP we make the change of variable
\begin{align*}
  \tilde{y}(\xib,t) &:= y(\xib,t) -y_d(\xib,t),&
  \tilde{x}(\xib)&:= x(\xib)-y_{d}(\xib,0).
\end{align*}
\(\tilde{y}\) minimizes the modified cost functional
\begin{align}
	\label{cost_functional}
	\tilde{\mathcal{J}}_{\tilde{x}}(\tilde{y}, u)
	&:= \frac{1}{2} \int_0^{\infty} (\|\tilde{y}(\cdot,\eta)\|_{L^{2}(\omg)}^{2}  + \gamma \|u(\eta)q(\cdot)\|_{L^{2}(\omg)}^2) e^{-\lambda \eta} d\eta
\end{align}
subject to the modified state equation
\begin{equation}
	\label{eq:test1_change}
	\left\{
	\begin{array}{rll}
		\partial_t \tilde{y}(\xib,t) + (-\Delta)^{\fracOrder}\tilde{y}(\xib,t) &= b(\xib,t) + u(t)q(\xib) \; -& (\xib,t) \in \omg \times (0,\infty){,} \\
		& \qquad -\partial_t y_d(\xib,t) - (-\Delta)^{\fracOrder}y_d(\xib,t) & \\
		\tilde{y}(\xib,t)&=0 &(\xib,t)\in \omg^{c} \times (0,\infty), \\
		\tilde{y}(\xib,0)&=\tilde{x}(\xib) &\xib\in \omg .
	\end{array}
	\right.
\end{equation}
Equation \eqref{eq:test1_change} is equivalent to \eqref{eq:test1}. The great benefit is that the modified cost functional has only the norm of $\tilde{y}$. With this change of variables we can solve the problem using the HJB approach and then recover the solution $y$.

The finite element spatial discretization of equation \eqref{eq:test1_change} leads us (under the same abuse of notation as in Section~\ref{sec:problemSetting}) to the system of ODEs
\begin{equation}
	\label{eq:test1_discrete}
	\left\{
	\begin{array}{rll}
		M \dot{\tilde{y}}(t) &= -A \tilde{y}(t) +  u(t) Q  + B(t) - M\dot{y_{d}}(t) - Ay_{d}(t) & t \in (0,\infty)\\
		\tilde{y}(0) &= \tilde{x},
	\end{array}
	\right.
\end{equation}
with $A,M \in \mathbb{R}^{d\times d}$ and $d$-dimensional vectors $\dot{\tilde{y}}, \tilde{y}, Q$ and $B$ as described in \eqref{fin_element}. Here, $y_d$ and $\dot{y_{d}}$ are also $d$-dimensional vectors with construction analogue to the other terms. We solve the closed loop problem with \eqref{eq:test1_discrete} as the state equation in order to approximate \eqref{cost_funct1}, and cost functional \eqref{cost_functional}. 

In this example, we consider $d=63$ and $d=127$. We compare the solution obtained by the DP approach, $y_{HJB}$, with the exact solution $y^*$ and with the solution obtained considering the exact optimal control plugged into the discrete system \eqref{eq:test1_discrete}, $y(u^{*})$.
To approximate the value function and obtain the control in feedback form, we consider $s=0.75, \gamma = 0.01, \lambda =0.5$ and $T_0=3$. To generate the scattered mesh the control space $U=[0,1]$ is discretized with $7$ equidistributed controls and the time step used is $\overline{\Delta t} = 0.0125$. The problem \eqref{eq:test1_discrete} is solved up to $T=4$ and the discrete trajectories are collected. The grid is formed with a total of $2248$ nodes for both $d=63$ and $d=127$. The separation distances are $0.0248$ and $0.0350$, respectively.

We use the grid generated by the dynamics to run our algorithm using the parameter space $\mathcal{P}=[0.1, 0.3]$ discretized with step size $0.02$. The residual is minimized with the parameter equal to 0.2. We consider $U$ discretized in $21$ points and $\Delta t = 0.01$. We use these data and the grid $X$ to approximate the value function and compute the feedback control. To obtain the feedback control and optimal trajectories we discretize the control set $U$ in $1681$ points and solve \eqref{eq:test1_discrete} up to the final time $T=4$. We invert the change of variables used in equation \eqref{eq:test1_change} and we compare the solutions obtained via DP approach, denoted by $y_{HJB}$, with the exact solution $y^*$ (as described in equation \eqref{eq:test1_optimal}) and with the solution obtained using the optimal control $u^{*}$ inserted on equation  \eqref{eq:test1_discrete}, denoted by $y(u^{*})$.


%
%

\begin{figure}[htbp]
	\centering 	
	
	\includegraphics[scale = 0.26]{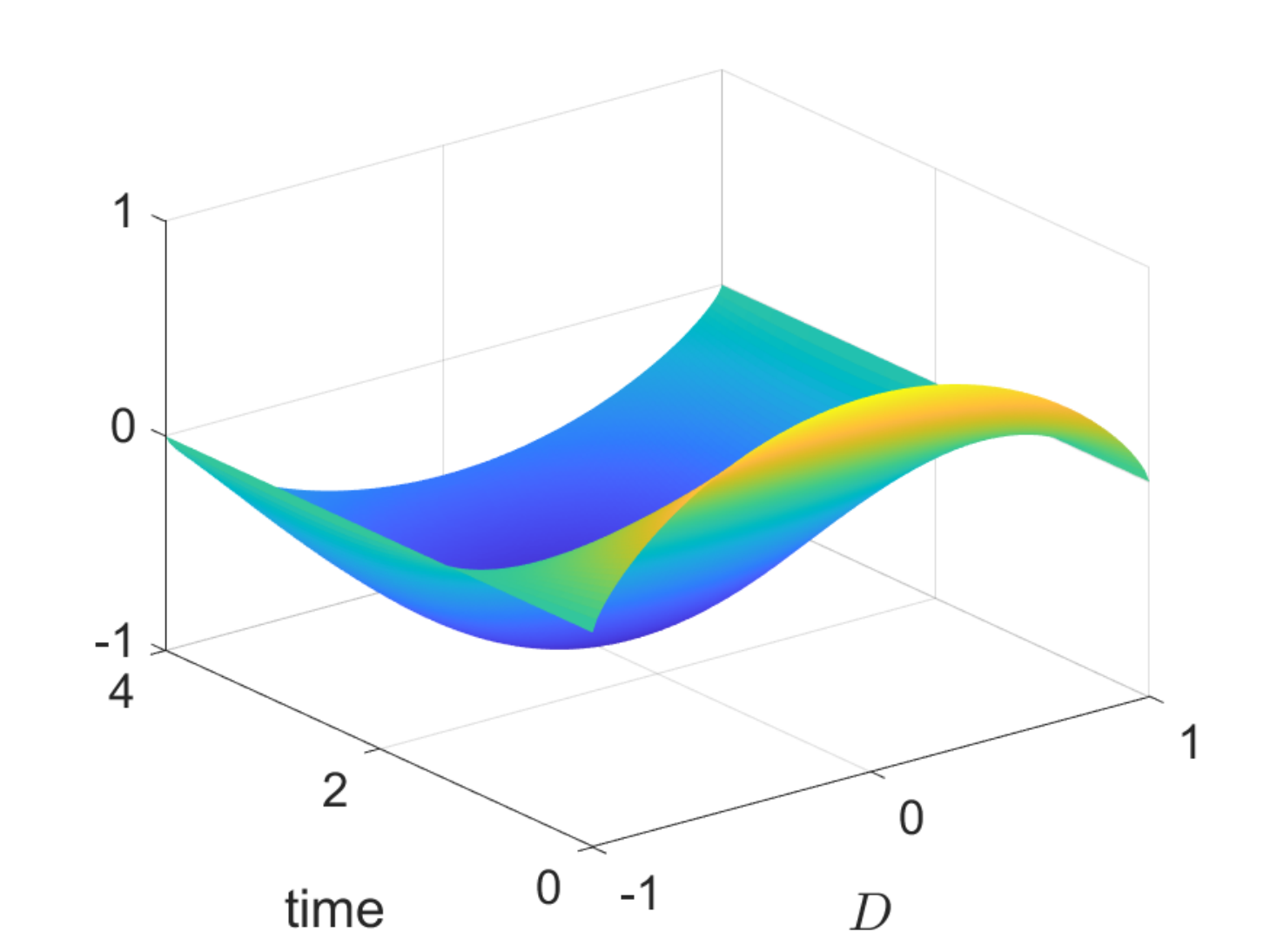}
	\includegraphics[scale = 0.26]{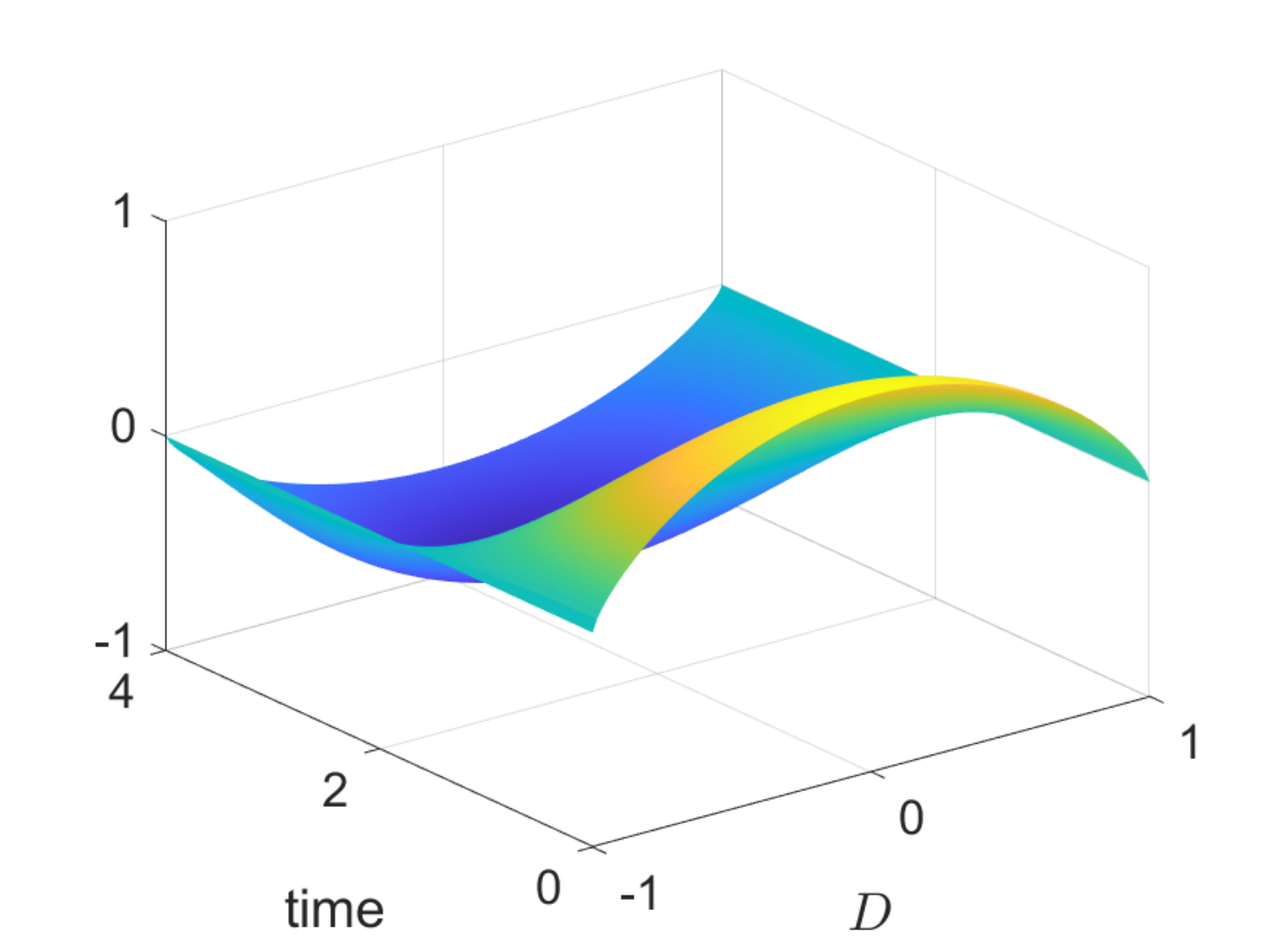}
	\includegraphics[scale = 0.26]{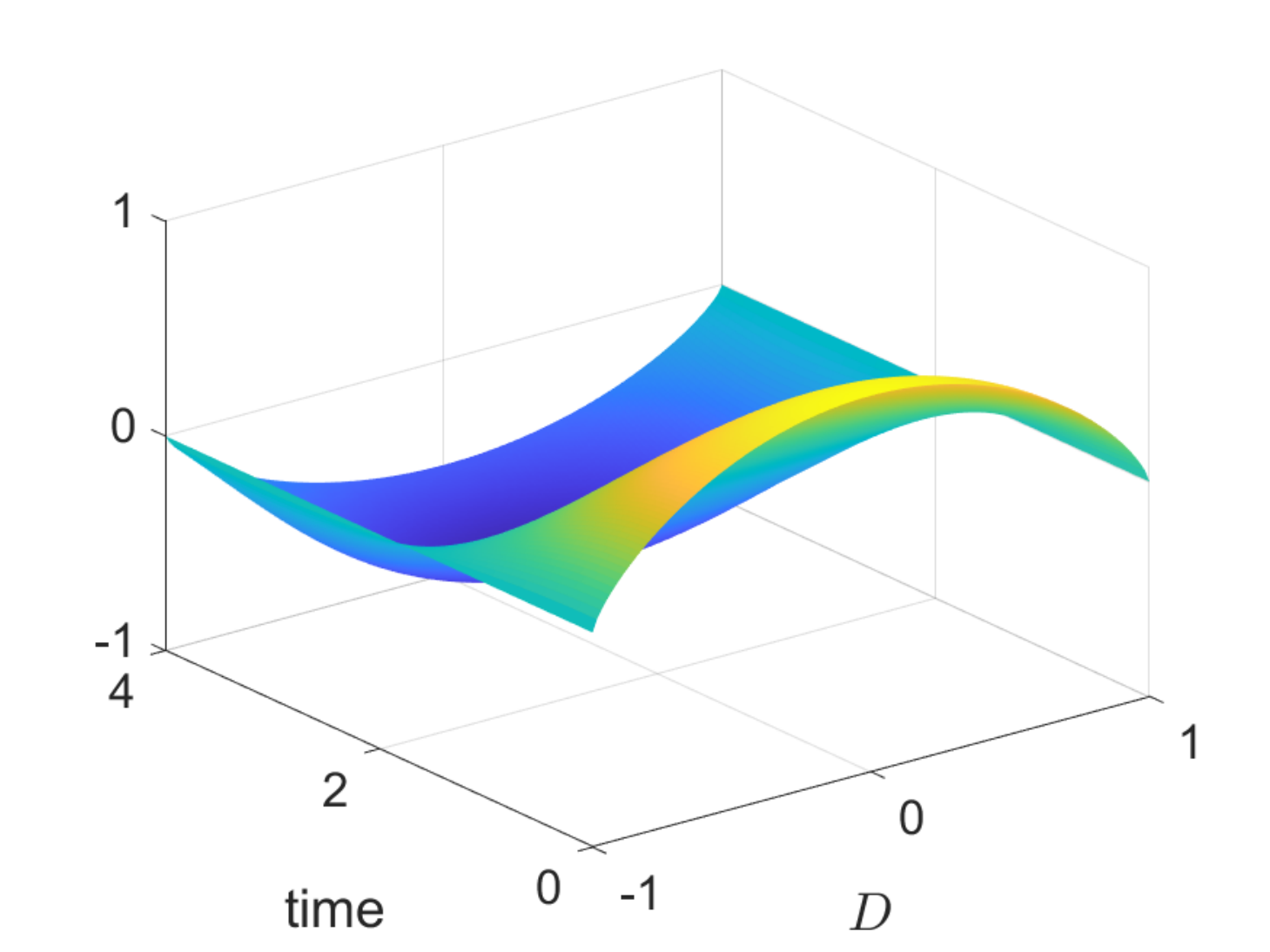}
	
	\caption{Test 1.
          Left: Uncontrolled solution (i.e. for \(u=0\)).
          Middle: Analytical solution.
          Right: HJB solution \(y_{HJB}\). Solutions in the case $\Delta t = 0.0125$ and $d=127$. }
	\label{fig1_test1}
\end{figure}
\begin{figure}[htbp]
	\centering 	
	\includegraphics[scale = 0.26]{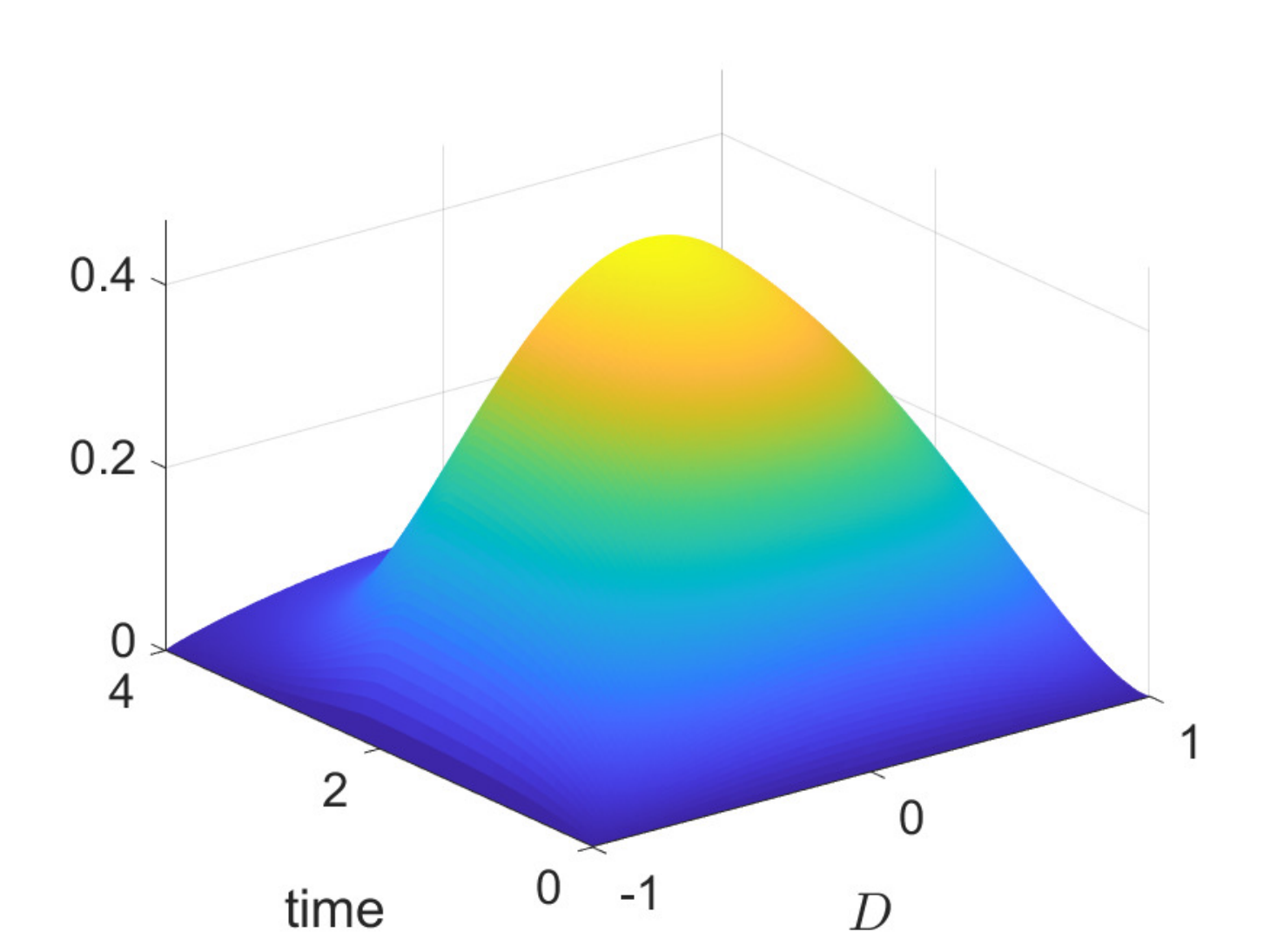}
	\includegraphics[scale = 0.26]{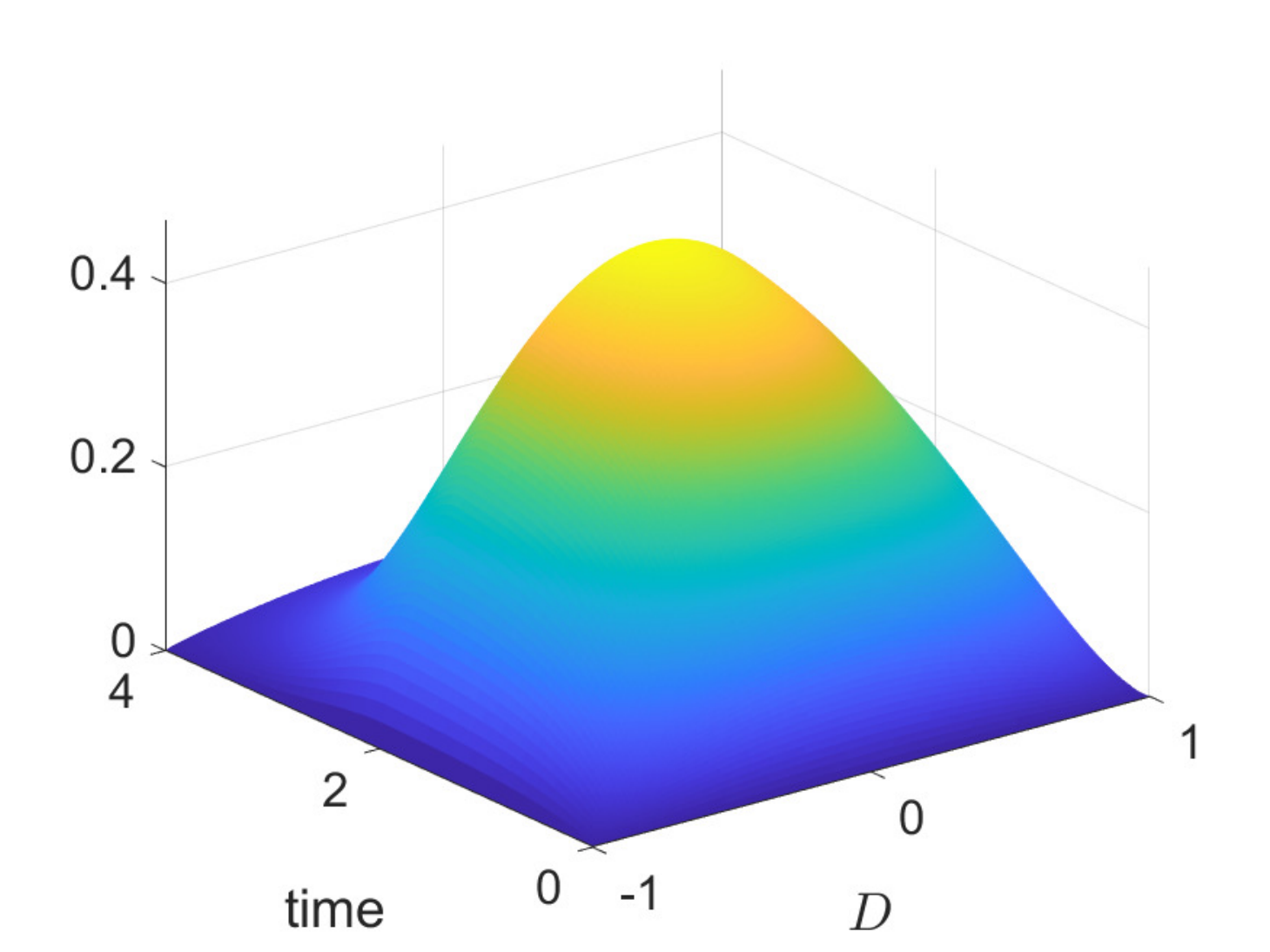}
	\includegraphics[scale = 0.26]{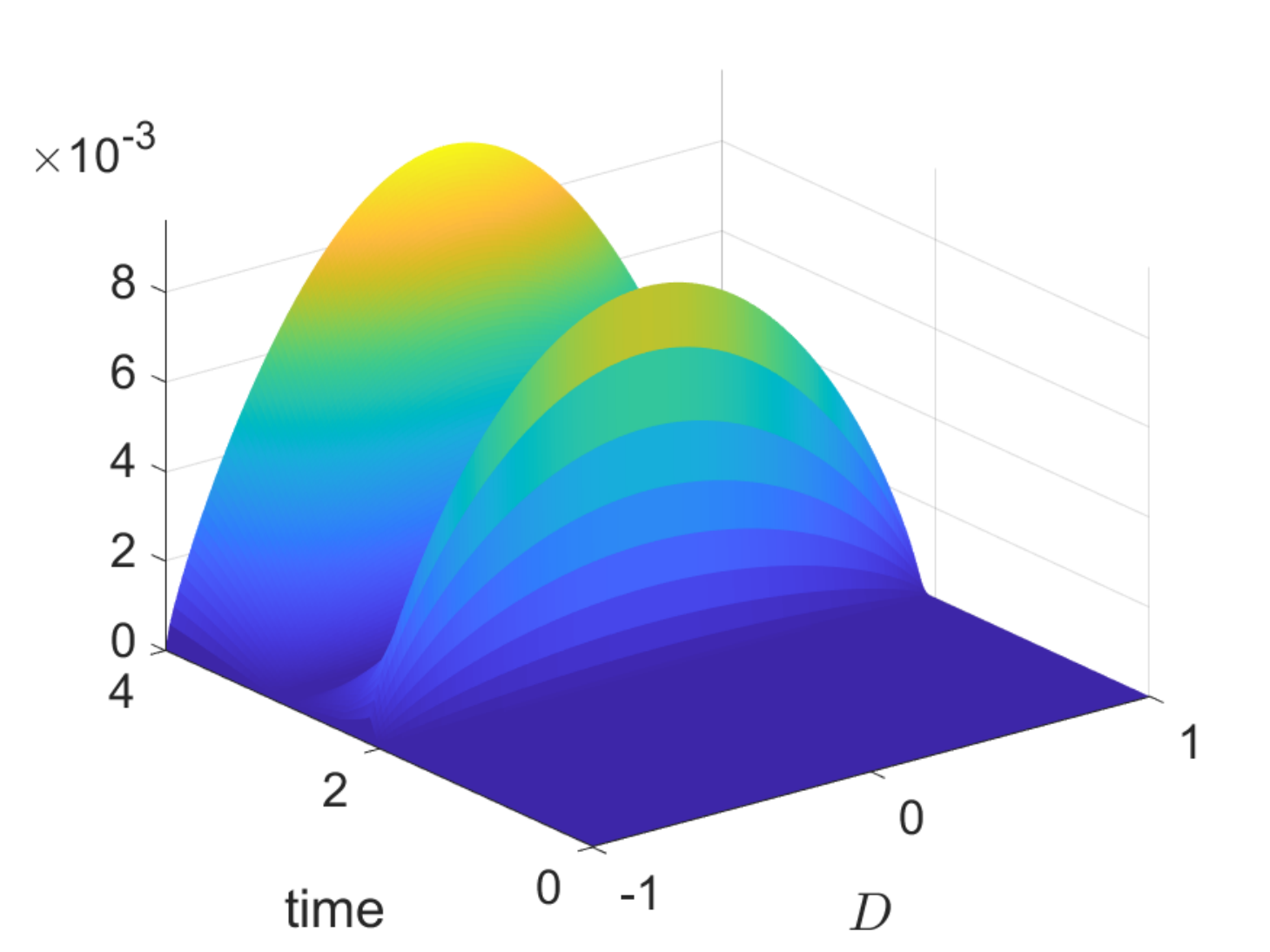}
	\caption{Test 1.
          Left: Absolute difference between the HJB solution and the uncontrolled solution.
          Center:  Absolute difference between the analytical solution  and the uncontrolled solution.
          Right: Absolute difference between the HJB and analytical solutions. }
	\label{abs_diff}
\end{figure}

Figure \ref{fig1_test1} presents the solutions for $d= 127$ and $\Delta t = 0.0125$. In the left panel we report the uncontrolled solution (i.e. \(u=0\)) of system \eqref{eq:test1_discrete}; in the middle panel we report the analytical solution. In the right panel we show the controlled solution via the HJB equation. We note the difference between the shape of the uncontrolled solution and the other solutions. This difference is clear in Figure \ref{abs_diff} where we display the absolute difference between the uncontrolled solution and the HJB solution and the difference between the uncontrolled solution and the analytical solution.  Also, the similarity between the shape of the analytic solution and the solution obtained by feedback control becomes evident when we look at their absolute difference in the right panel of Figure \ref{abs_diff} (note that the order of the $z-$axis is $O(10^{-3})$).

\begin{figure}[htbp]
	\centering 	
	\includegraphics[scale = 0.39]{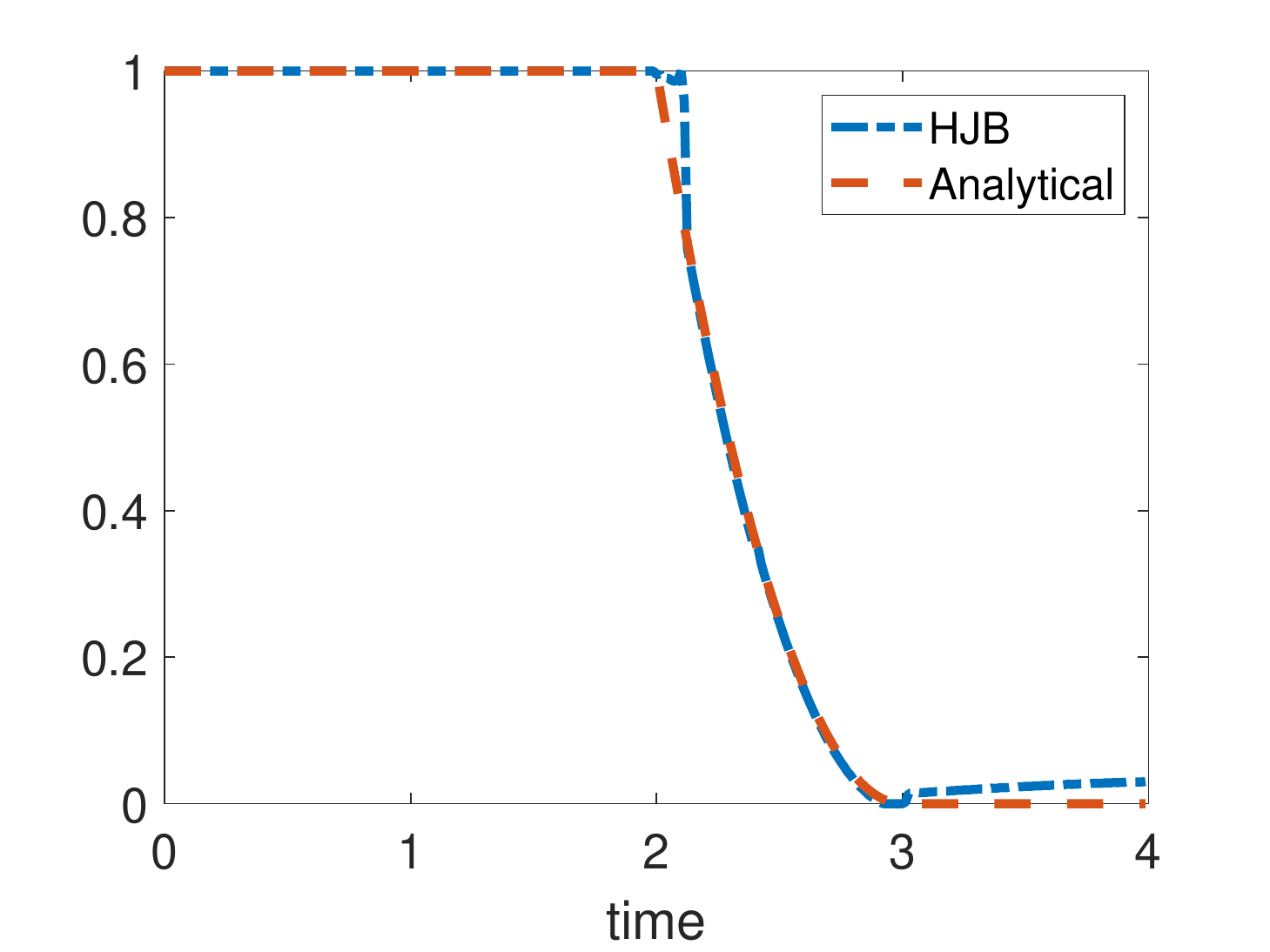}
	\includegraphics[scale = 0.39]{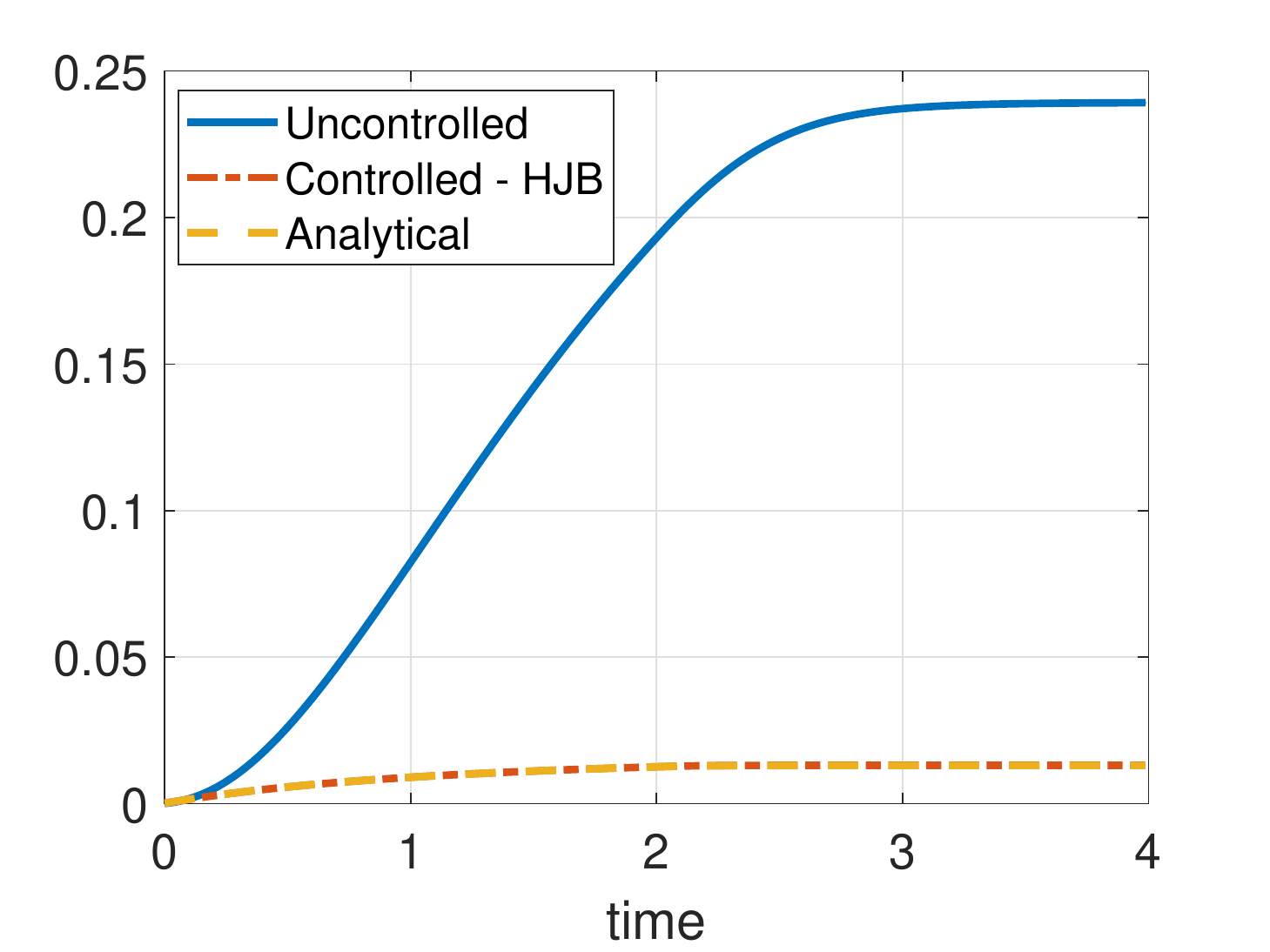}
 	\caption{Test 1. Cost functional and absolute differences in the case $\Delta t = 0.0125$ and  $d=127$.
          Left: optimal control.
          Right: Cost functional \eqref{cost_funct1} calculated for uncontrolled, analytical and controlled by HJB approach solutions.}
 	
	\label{fig2_test1}
\end{figure}

In the left panel of Figure \ref{fig2_test1}, we compare the control obtained by the HJB approach and the analytical control $u^{*}$. We observe that the controls are very close up to $T_0$. Then, there is a slight difference which is not significant enough to affect the value of the cost functional.
Indeed, in the right panel of Figure \ref{fig2_test1} we show the evaluation of the cost functional for the uncontrolled solution, the analytical and the HJB driven approach. Here, we observe that the cost functional of the two controlled solutions matches perfectly. This highlights the quality of the proposed approximation.





To further validate our approach, we study the error of the proposed approximation with respect to the control of the nonlocal equation.
The convergence analysis is performed with the norm  $\|\cdot\|_{L_{\nu}^{2}(0,T;\omg)}$, which we write as $\|\cdot\|_{L_{\nu}^{2}}$ for ease of notation.
In Table \ref{test1_table63} we show the values of $\|y_{HJB}-y^{*}\|_{L_{\nu}^{2}}$ in the second column, $\|y_{HJB}-y(u^{*})\|_{L_{\nu}^{2}}$ in the fourth column and $\|y^{*}-y(u^{*})\|_{L_{\nu}^{2}}$ in the sixth. Here we considered distinct values of $\Delta t$ to solve the discrete problem \eqref{eq:test1_discrete}. We note that the errors decay with the reduction in $\Delta t$, as expected, since a more refined time discretization implies better approximations.
Table~\ref{test1_table63} shows that the spatial discretization error of the state equation is dominated by temporal discretization error (column 6).
In turn, the temporal discretization error of the state equation is dominated by the discretization error of the HJB solution.

\begin{table}[htbp]
	\centering
	\begin{tabular}{ccccccc}
		\hline
		$\Delta t$	& $\|y_{HJB}-y^{*}\|_{L_{\nu}^{2}}$ & rate & $\|y_{HJB}-y(u^{*})\|_{L_{\nu}^{2}}$ & rate & $\|y^{*}-y(u^{*})\|_{L_{\nu}^{2}}$ & rate \\
		\hline
		$0.05$	& $0.0450$ & & $0.0351$  & & $0.0470$ & \\
		$0.025$	& $0.0340 $ & $0.40$ & $0.0225$  & $0.64$ &$0.0233$ & $1$\\
		$0.0125$ & $0.0260$ & $0.38$ & $0.0197$ & $0.2$ &$0.0115$ & $1$\\
		\hline
	\end{tabular}  
	\caption{Test 1. Data of simulation with dimension = $63$, $h=0.0248$ and number of nodes = $2248$.}	
	\label{test1_table63}
\end{table}

Table \ref{test1_table127} presents results for $d=127$. The decay in the norm of differences also happens with a reduction of the time step. The comparison between values in Table \ref{test1_table127} with values in Table \ref{test1_table63} suggests that the refined spatial discretization implies a more accurate approximation, since the values in the second and fourth columns are smaller than their counterpart in Table \ref{test1_table63}. These values in Table \ref{test1_table127} are almost the half of the respective values in Table \ref{test1_table63}. The entries in the sixth column of Table \ref{test1_table63} are very similar to those in the same column of Table \ref{test1_table127}, indicating that the time discretization error has a greater influence than the spatial error.

\begin{table}[htbp]
	\centering
	\begin{tabular}{ccccccc}
		\hline
		$\Delta t$	& $\|y_{HJB}-y^{*}\|_{L_{\nu}^{2}}$ & rate & $\|y_{HJB}-y(u^{*})\|_{L_{\nu}^{2}}$ & rate & $\|y^{*}-y(u^{*})\|_{L_{\nu}^{2}}$ & rate\\
		
		\hline
		$0.05$	& $0.0423$ & & $0.0220$  & & $0.0473$& \\
		$0.025$	& $0.0281$ & $0.59$ & $0.0173$ & $0.34$ &$0.0236$ & $1$\\
		$0.0125$ & $0.0149$ & $0.9$ & $0.0072$ & $1.2$ & $0.0117$ & $1$ \\
		
		\hline
	\end{tabular}  
	\caption{Test 1. Data of simulation with dimension $d=127$, $h=0.0350$ and number of nodes = $2248$.}	
	\label{test1_table127}
\end{table}



		
%
\paragraph{Simulation with noise}
One of the benefits of computing the control in closed form is the ability to react under disturbances of the system and changes in the initial condition. In this example, we control the system with initial condition $x=2q$ and disturb the system at each time instance by adding a noise vector of independent, identically distributed Gaussian variables with zero mean and standard deviation $0.025$. 

\begin{figure}[H]
	\centering 	
	\includegraphics[scale = 0.35]{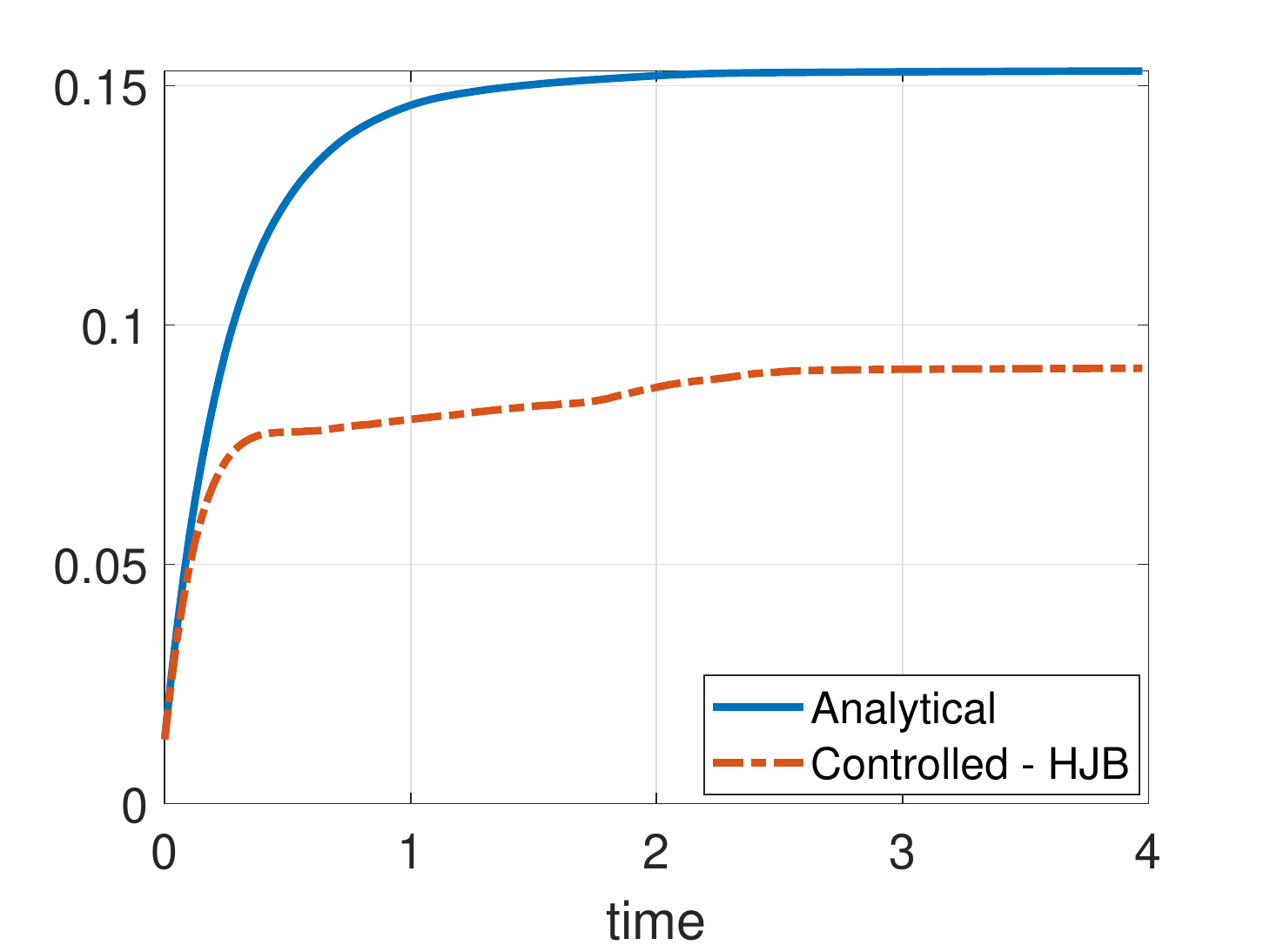}
		\includegraphics[scale = 0.35]{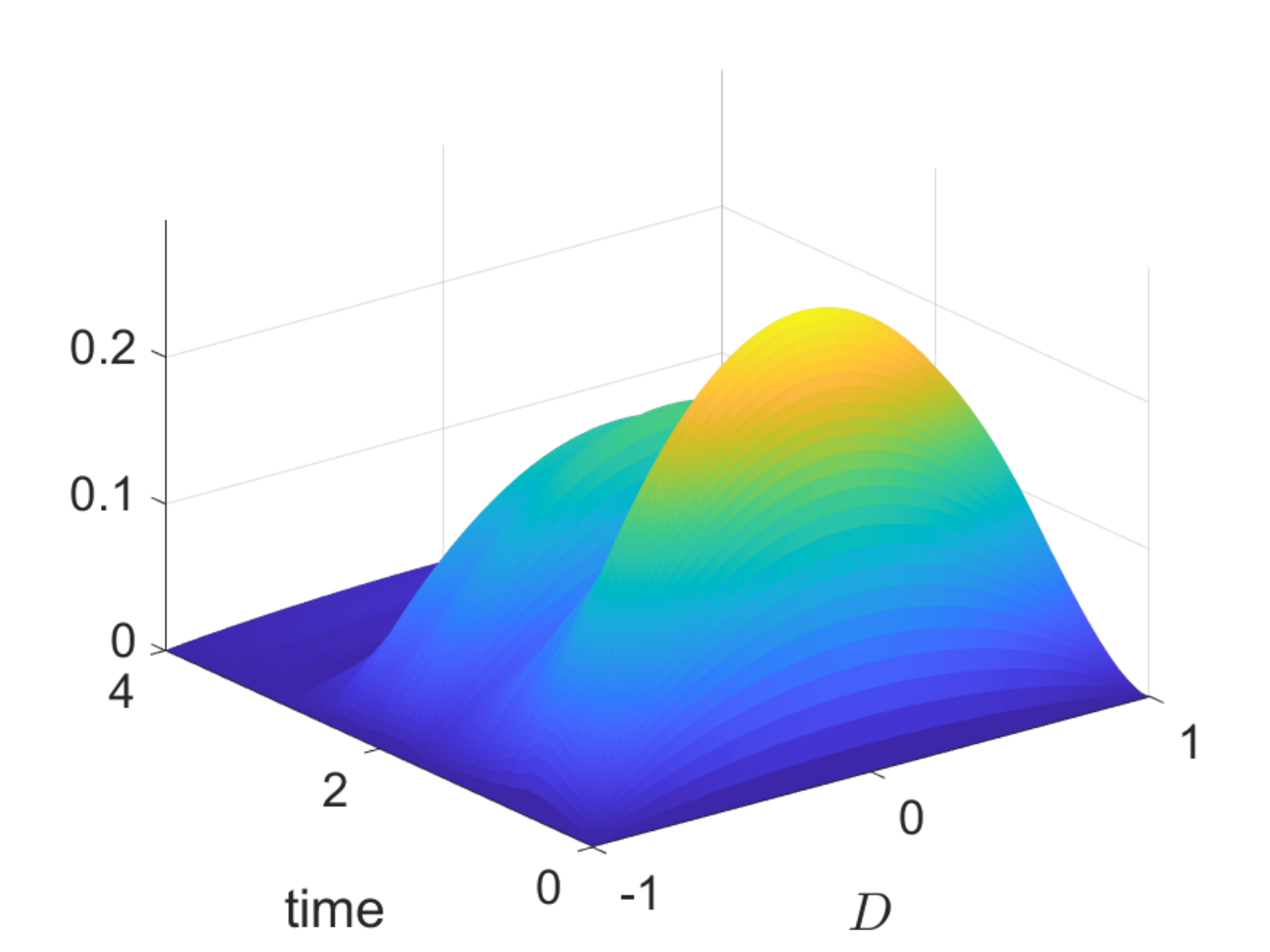}\\
	\includegraphics[scale = 0.35]{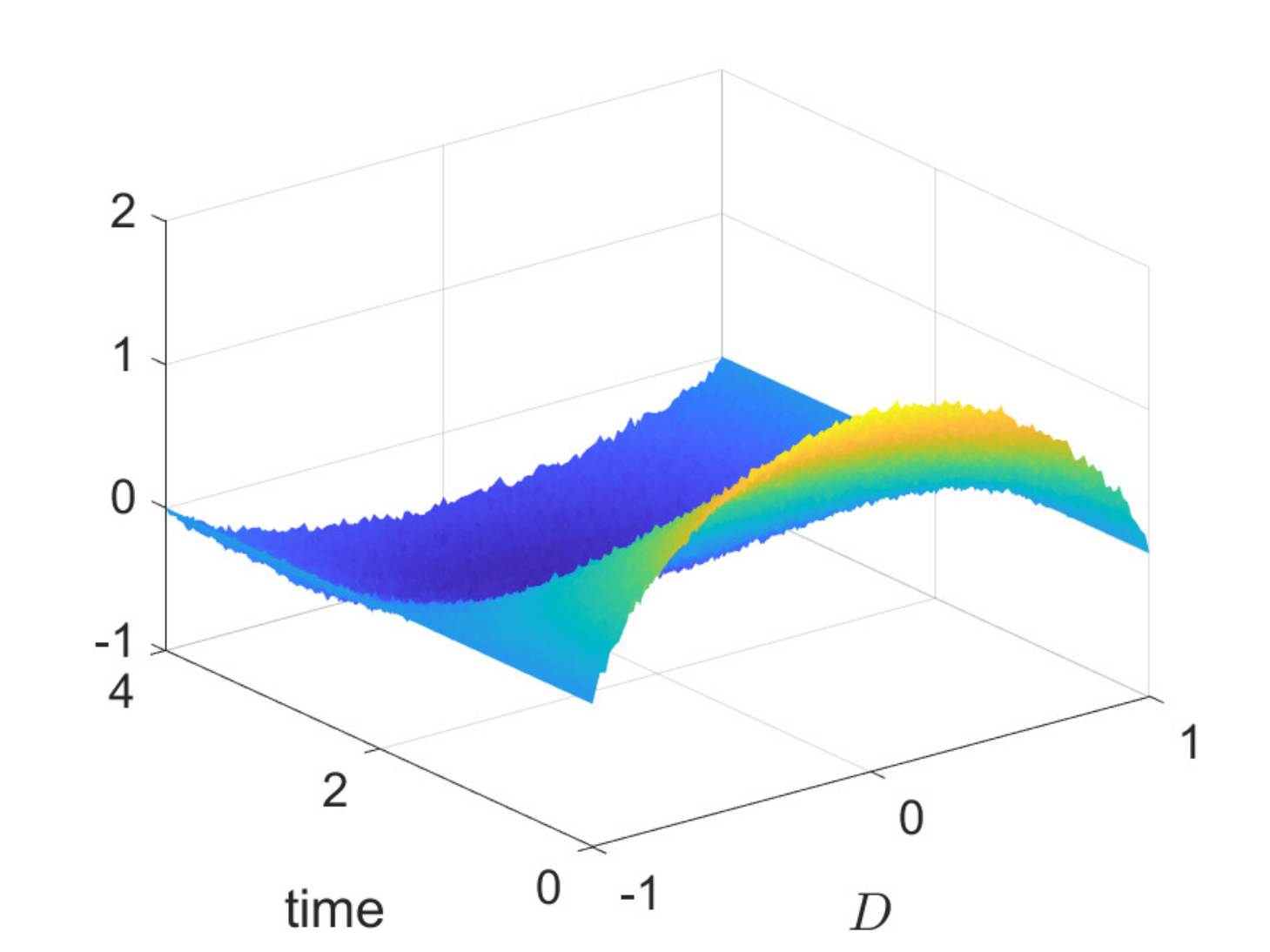}	\includegraphics[scale = 0.35]{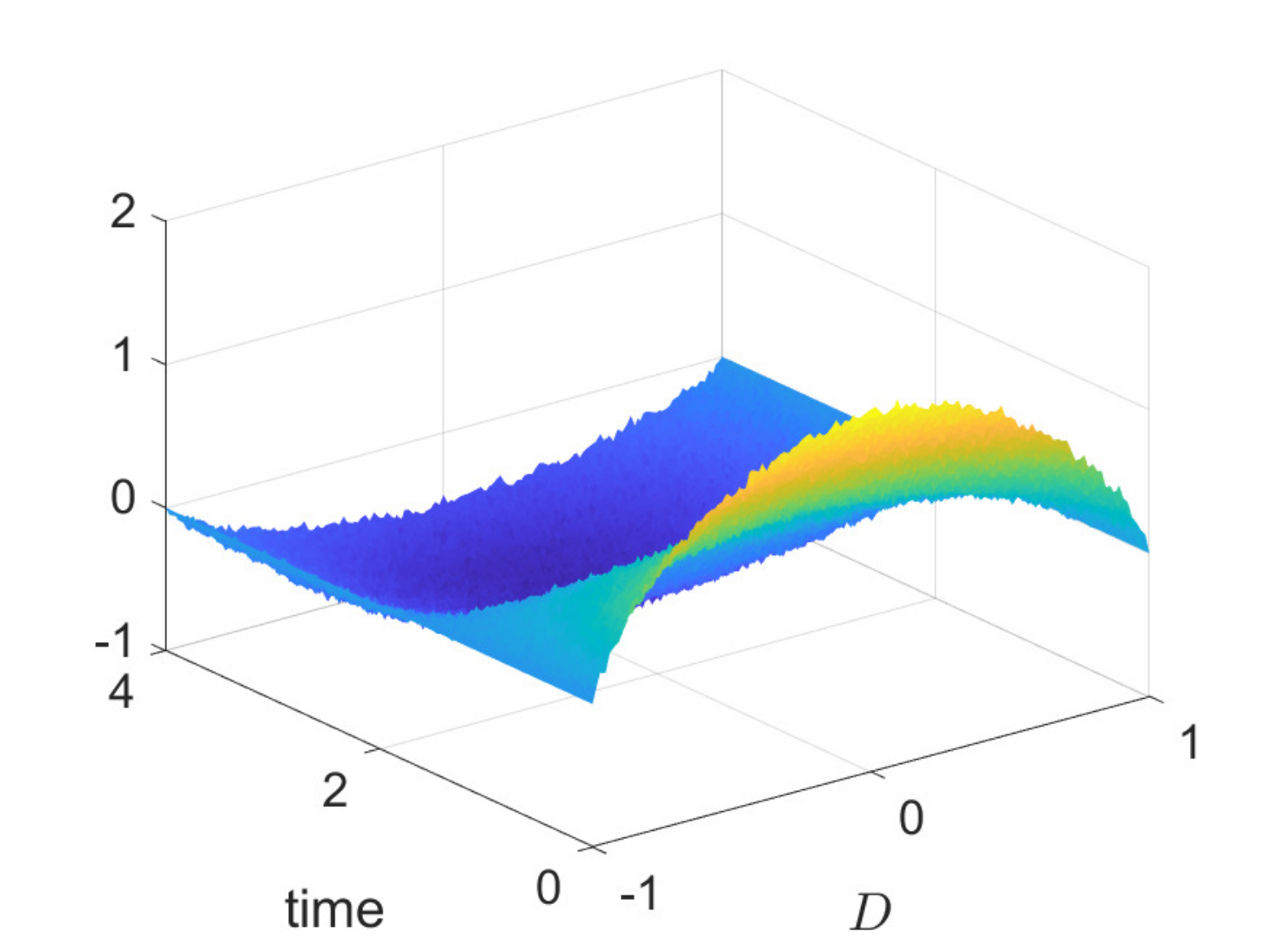}
		\caption{Test 1. Top: evaluation of the cost functionals with noise term (left), absolute difference between HJB and controlled solution with analytical control (right). Bottom:  Solution of the disturbed system with analytical control. (left), HJB approach solution with noise (right).} 
	\label{fig_noise}
\end{figure}

Then, with the same value function obtained in the previous test and $d=127$, we compute the feedback control with $\Delta t = 0.025$ and \(U\) discretized with $1681$ nodes. We compare the solution using the analytical control plugged into the disturbed system and the solution using the HJB approach which adjusts the control based on disturbance of the system. We can see, in the top left panel of Figure \ref{fig_noise}, the evaluation of the cost functional for both methods. It is clear that the feedback control leads to a lower cost functional than the pre-determined control. For completeness, we show the controlled solutions in the bottom panels of Figure \ref{fig_noise} where we report the analytical control on the left and the HJB control on the right.
In the top right panel of Figure \ref{fig_noise} we report the absolute difference between both solutions.

This test shows the robustness of our method under disturbances and motivates the use of this approach as opposed to open-loop control.

\subsection{Test 2: Linear example with control acting on target domain}
In this example we minimize the cost functional
\begin{align}
	\mathcal{J}_{x}(y, u)
	&:= \frac{1}{2} \int_0^{\infty} (\|y(\cdot,\eta)\|_{L^{2}(\mathcal{T})}^{2}  + \gamma \|q(\cdot)u(\eta)\|_{L^{2}(\omg)}^2) e^{-\lambda \eta} d\eta;
\end{align}
here, the main difference with respect to the cost functional \eqref{cost_functional}, is that the norm of $y$ is restricted to the target $\mathcal{T}=[-1/2,1/2]$ with $\mathcal{T} \subset \omg$ and $\omg=(-1,1)$. This functional is minimized subject to
\begin{equation}  	
	\label{eq:test2}
	\left\{
	\begin{array}{rll}
		\partial_t y(\xib,t) + (-\Delta)^{\fracOrder}y(\xib,t) &= b(\xib,t) + u_{1}(t)q_{1}(\xib) + u_{2}(t)q_{2}(\xib) &\text{for all }\xib\in \omg\times (0,T){,}  \\
		y(\xib,t)&=0 &\text{for all }\xib\in \omg^{c} \times (0,T) , \\
		y(\xib,0)&=x(\xib) &\text{for all }\xib\in \omg,
	\end{array}
	\right.
\end{equation}
where \begin{equation}
	b(\xib,t) = (1-\cos(t))\chi_{[-1,-3/4]}(\xib)
\end{equation} 
and 
\begin{align}
	\label{control_test2}
	q_{1}(\xib)=\chi_{[-3/4, -1/2]}(\xib) ,& & q_{2}(\xib)=\chi_{[1/2, 3/4]}(\xib),
\end{align}
with $u_{i}(t)\in U=[-0.5,0]$, $i \in \{1,2\}$, i.e., the control is formed by two independent terms, each of them acting on a specific region of the domain.
The semi-discrete system associated with equation \eqref{eq:test2} is given by
\begin{equation}
	\left\{
	\begin{array}{rll}
		M \dot{y}(t) &= - A y(t) +  u_{1}(t)Q_{1} + u_{2}(t)Q_{2} + B(t),& t \in (0,\infty),\\
		y(0) &= 0.
	\end{array}
	\right.
	\label{eq:test2_discrete}
\end{equation}
The cost functional to be minimized is
\begin{align}
	\mathcal{J}_{x}(y, u)
	&= \frac{1}{2} \int_0^{\infty} (\|y(\eta)\|_{L^{2}(\mathcal{T})}^{2}  + \gamma \|u_{1}(\eta)q_{1}(\cdot) + u_{2}(\eta)q_{2}(\cdot)\|_{L^{2}(\omg)}^2) e^{-\lambda \eta} d\eta
\end{align}
with $y(t)$ solution of \eqref{eq:test2_discrete} and $u_{i}(t)$ defined as \eqref{control_test2}.


We set the parameters of the equation as $\gamma = 10^{-6}$ and $\lambda = 0.5$. To build the scattered mesh  we use $\overline{\Delta t} = 0.025$  and control pairs $(u_1,u_2) \in U\times U$, with $U$ discretized in $5$ equidistributed controls, resulting in a total of $25$ control pairs. We solve equation \eqref{eq:test2_discrete} up to time $T=6$ and collect $6026$ nodes to form a $63$-dimensional grid with separation distance $0.0180$.
We consider the parameter $\theta = 0.01$, time step $\Delta t =0.01$ and use the same discrete control set of the mesh; we compute the value function using the unstructured grid. In order to execute the feedback reconstruction we used $\Delta t = 0.025$ and  $1681$ control points.

\begin{figure}[htbp]
	\centering 	
	\includegraphics[scale = 0.26]{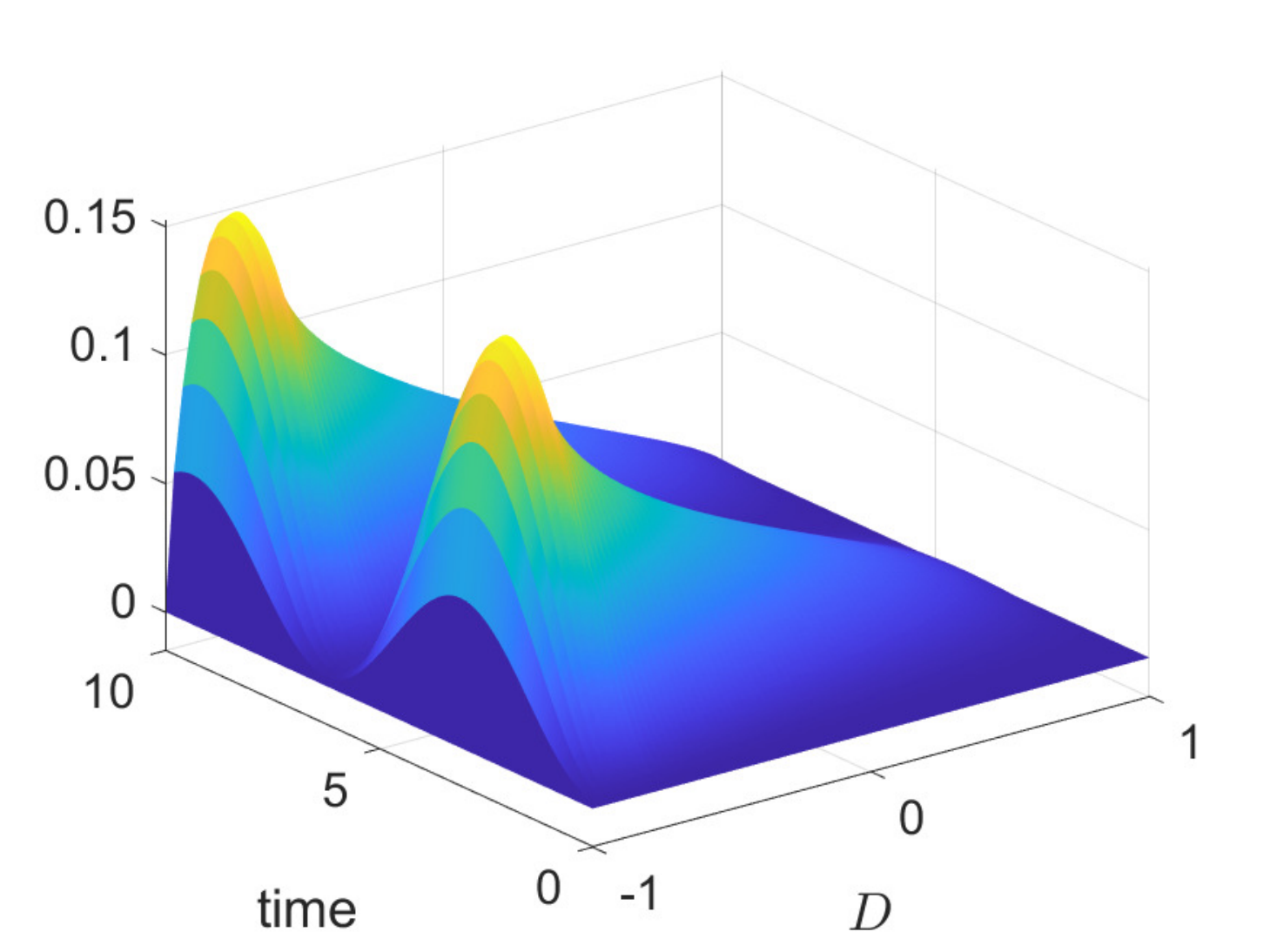}
	\includegraphics[scale = 0.26]{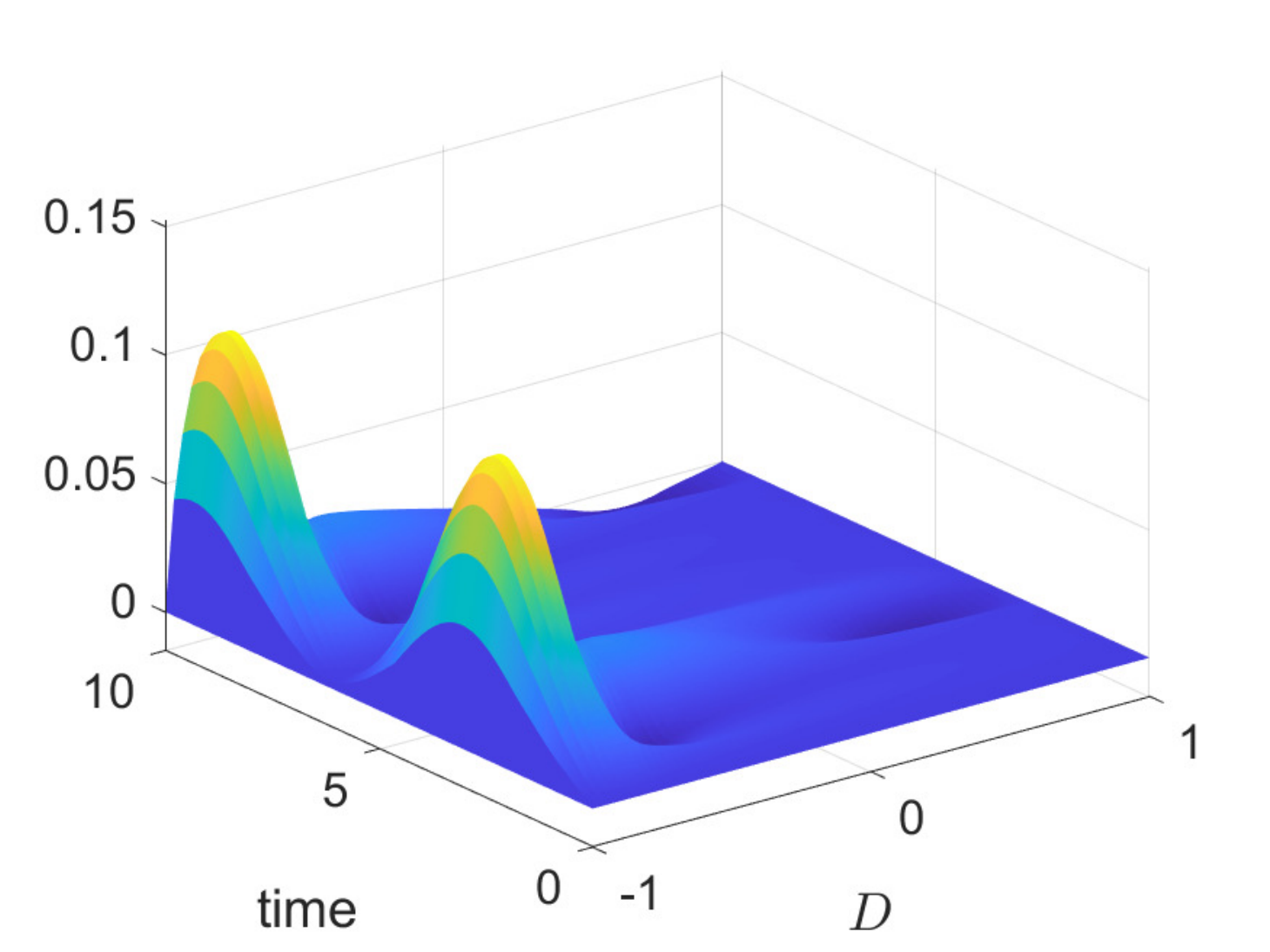}
	\includegraphics[scale = 0.26]{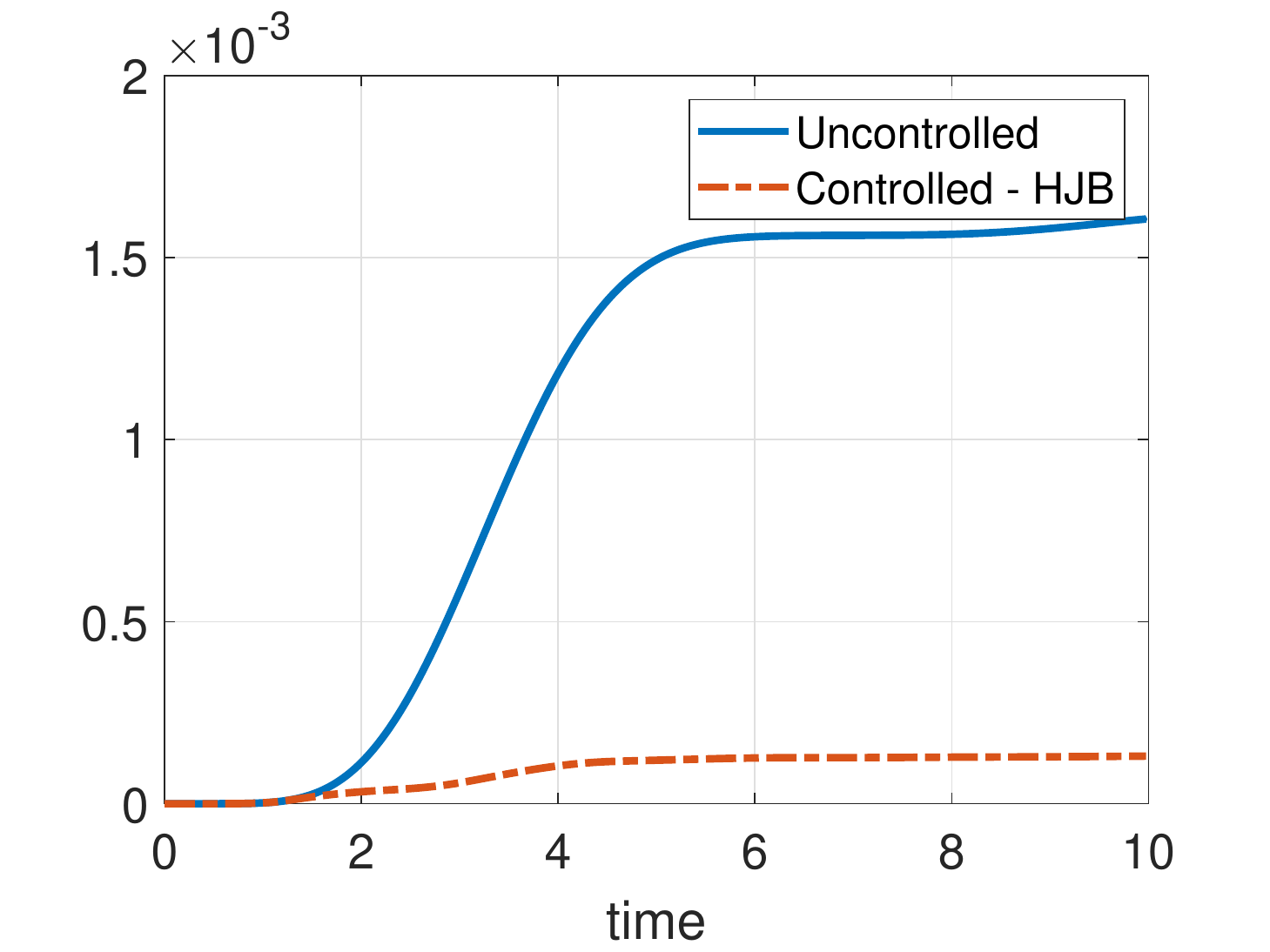}
	
	\caption{Test 2: Left: uncontrolled solution. Middle: controlled solution. Right: Cost Functionals. }
	\label{fig1_test2}
\end{figure}

In Figure \ref{fig1_test2}, we show the uncontrolled solution in the left panel and the controlled solution in the middle, both with final time $T=10$. The difference between the solutions becomes clear in the evaluation of the cost functional in the right panel of Figure \ref{fig1_test2}. The cost of the controlled solution is below the cost functional of the uncontrolled solution.

\begin{figure}[htbp]
	\centering 	
	\includegraphics[scale = 0.26]{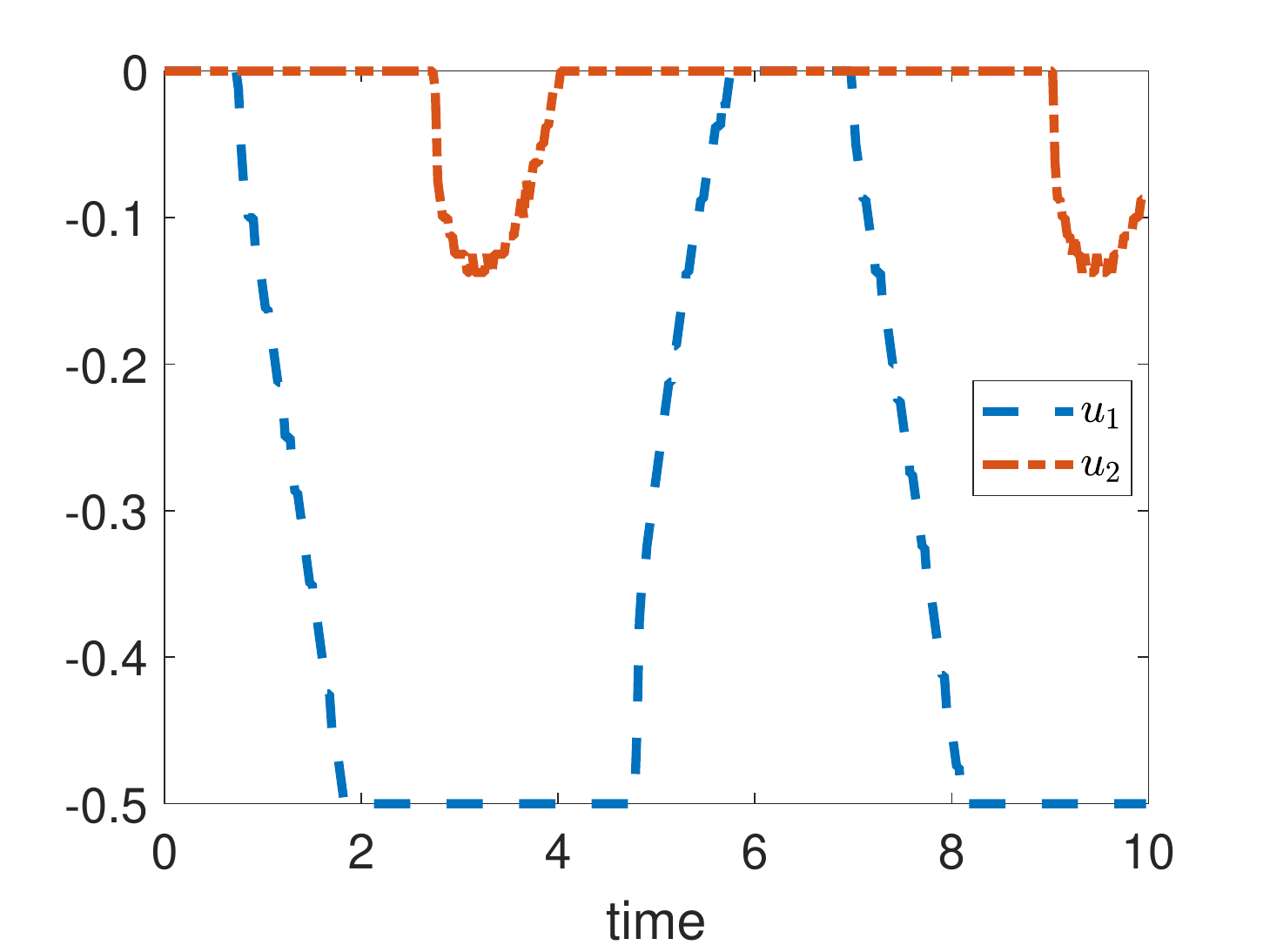}
	\includegraphics[scale = 0.26]{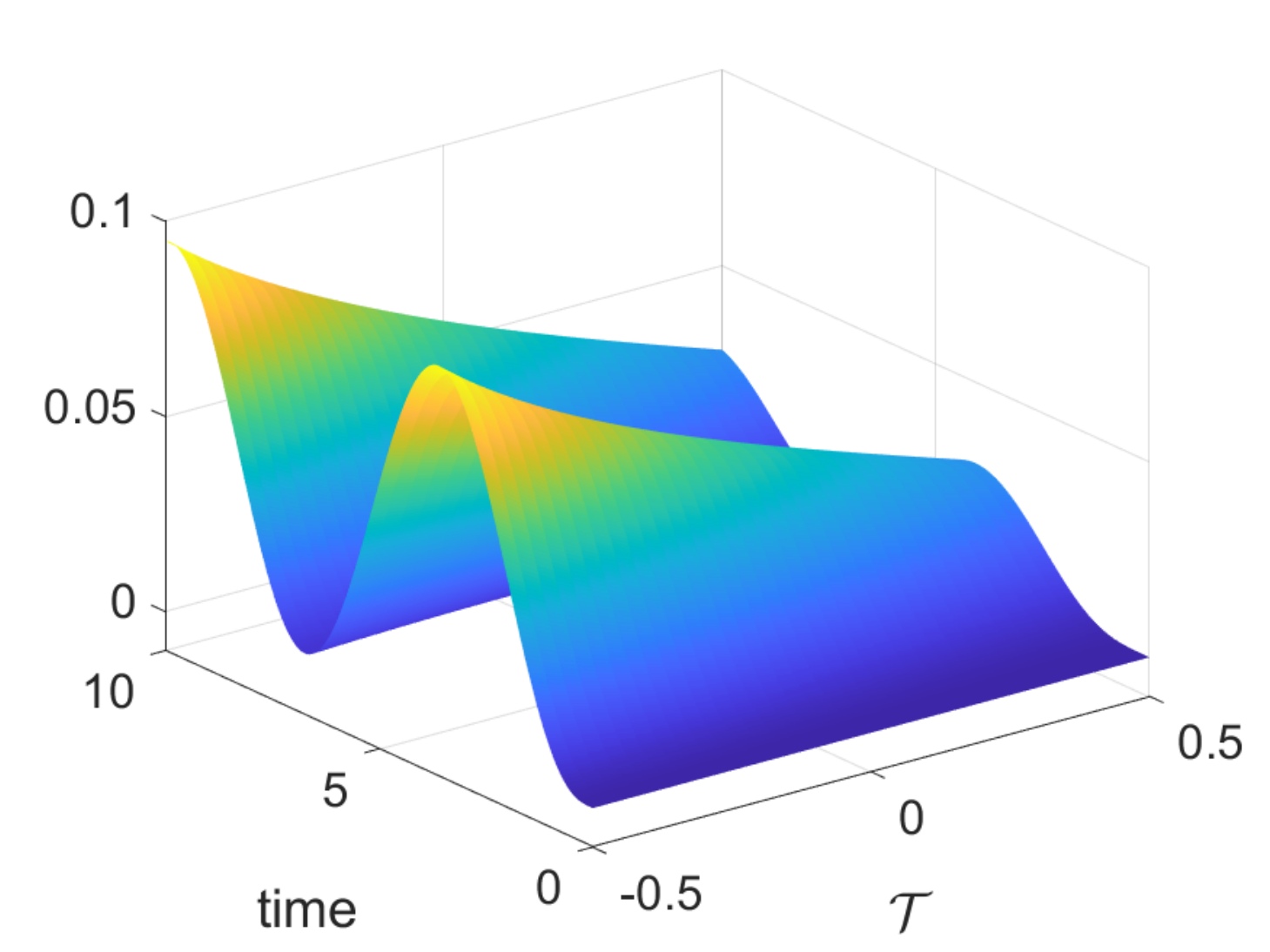}
	\includegraphics[scale = 0.26]{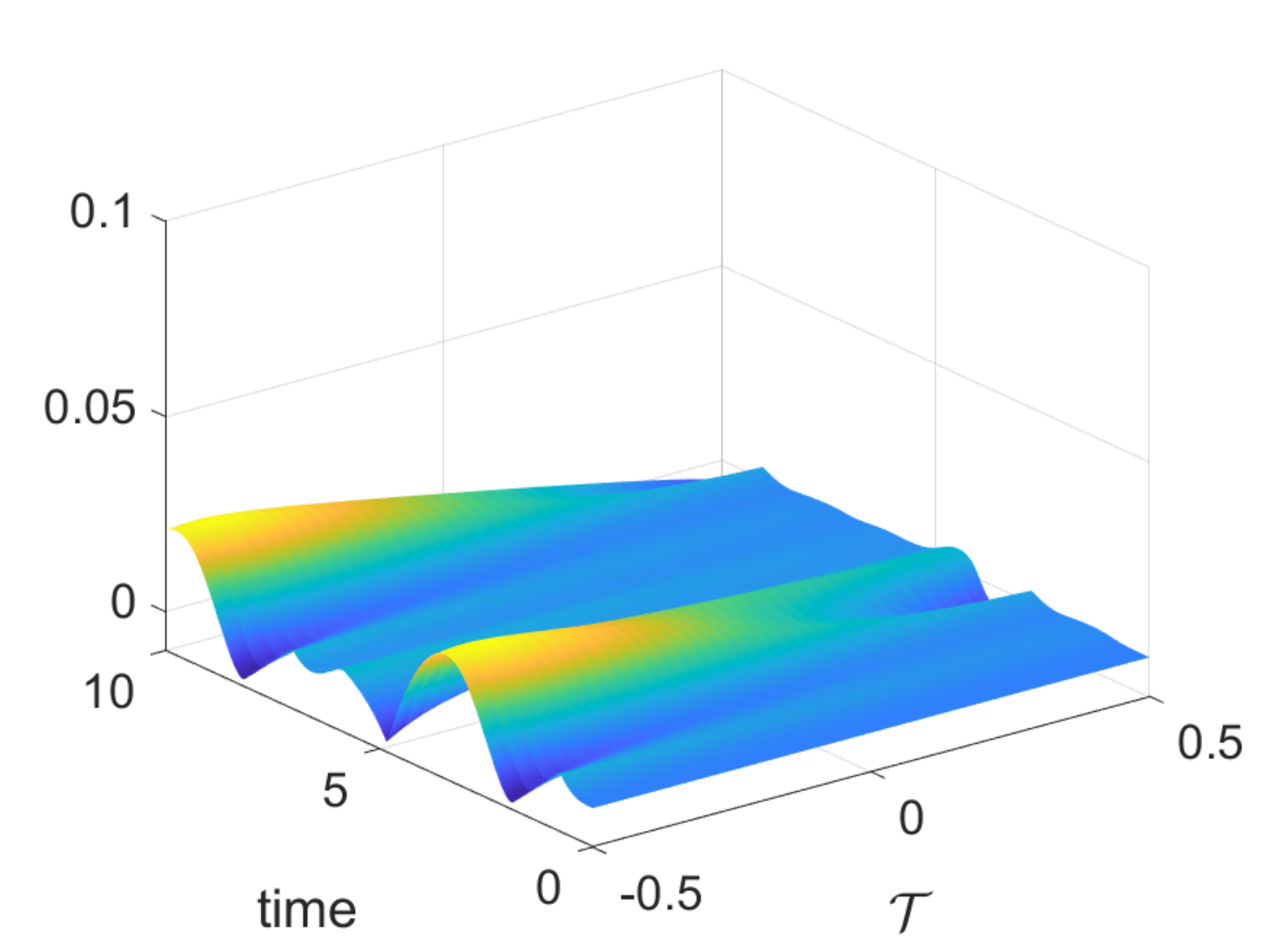}
	\caption{Test 2: Left: Controls $u_1$ and $u_2$. Center:  uncontrolled solution in a plot restricted to target $\mathcal{T}$. Right: controlled solution in a plot restricted to target $\mathcal{T}$.}
	\label{fig2_test2}
\end{figure}

In Figure \ref{fig2_test2}, we present the two control variables $u_1$ and $u_2$ obtained with our approach in the left panel. The control $u_1$ acts in an interval that is very close to the region where the force term is positive. Then, it is clear that $u_1$ is more active than $u_2$, as expected. Finally, a further zoom-in of the uncontrolled and controlled solution in the target region is reported in the middle and right panels, respectively.

\subsection{Test 3: A nonlinear example}

In this test we consider the following nonlinear state equation

\begin{equation}
	\left\{
	\begin{array}{rll}
		\partial_t y(\xib,t) + \alpha(-\Delta)^{\fracOrder}y(\xib,t) +  &=  F(y(\xib,t)) + u(t)q(\xib) &\text{for all }\xib\in \omg \times (0,\infty){,}  \\
		y(\xib,t)&=0 &\text{for all }\xib\in \omg^{c} \times (0,\infty) , \\
		y(\xib,0)&=x(\xib) &\text{for all }\xib\in \omg ,
	\end{array}
	\right.
\end{equation}
where the nonlinear term is given by $  F(y)= y^2(1-y)$ and the parameter $\alpha = 0.01$. The finite element spatial discretization results in the system
\begin{equation}
	\label{eq:test3_discrete}
	\left\{
	\begin{array}{rll}
		M \dot{y}(t) &= -\alpha A y(t)  + \mathbf{F}(y(t)) + u(t) Q,& t \in (0,\infty)\\
		y(0) &= x ,
	\end{array}
	\right.
\end{equation}
where the nonlinear term $\mathbf{F}: \mathbb{R}^d \rightarrow \mathbb{R}^d$ is $\mathbf{F}(y(t)) = y(t)^2 - y(t)^3$.
The initial condition used in this case is the same used in Section \ref{sec:example}, i.e. the solution of the Poisson problem.
The cost functional we minimize is given by
\begin{align}
	\mathcal{J}_{x}(y, u)
	:= \frac{1}{2} \int_0^{\infty} (\|y(\eta)\|_{L^{2}(\omg)}^{2}  + \gamma \|u(\eta)q(\cdot)\|_{L^{2}(\omg)}^2) e^{-\lambda \eta} d\eta,
\end{align}
subjected to the semi-discrete system \eqref{eq:test3_discrete}.
The parameters chosen in the cost are $\lambda = 0.5$ and $\gamma = 0.01$. To generate the scattered mesh we consider the control space $U=[-0.5,0]$ discretized in $11$ nodes, a temporal step $\overline{\Delta t} = 0.025$ and the same initial condition as for the first test. We collect the trajectory points generated solving the system \eqref{eq:test3_discrete} up to the final time $T=6$ for $d=63$. The scattered mesh has $2652$ nodes and the separation distance is $0.12$.

\begin{figure}[htbp]
	\centering 	
	\includegraphics[scale = 0.39]{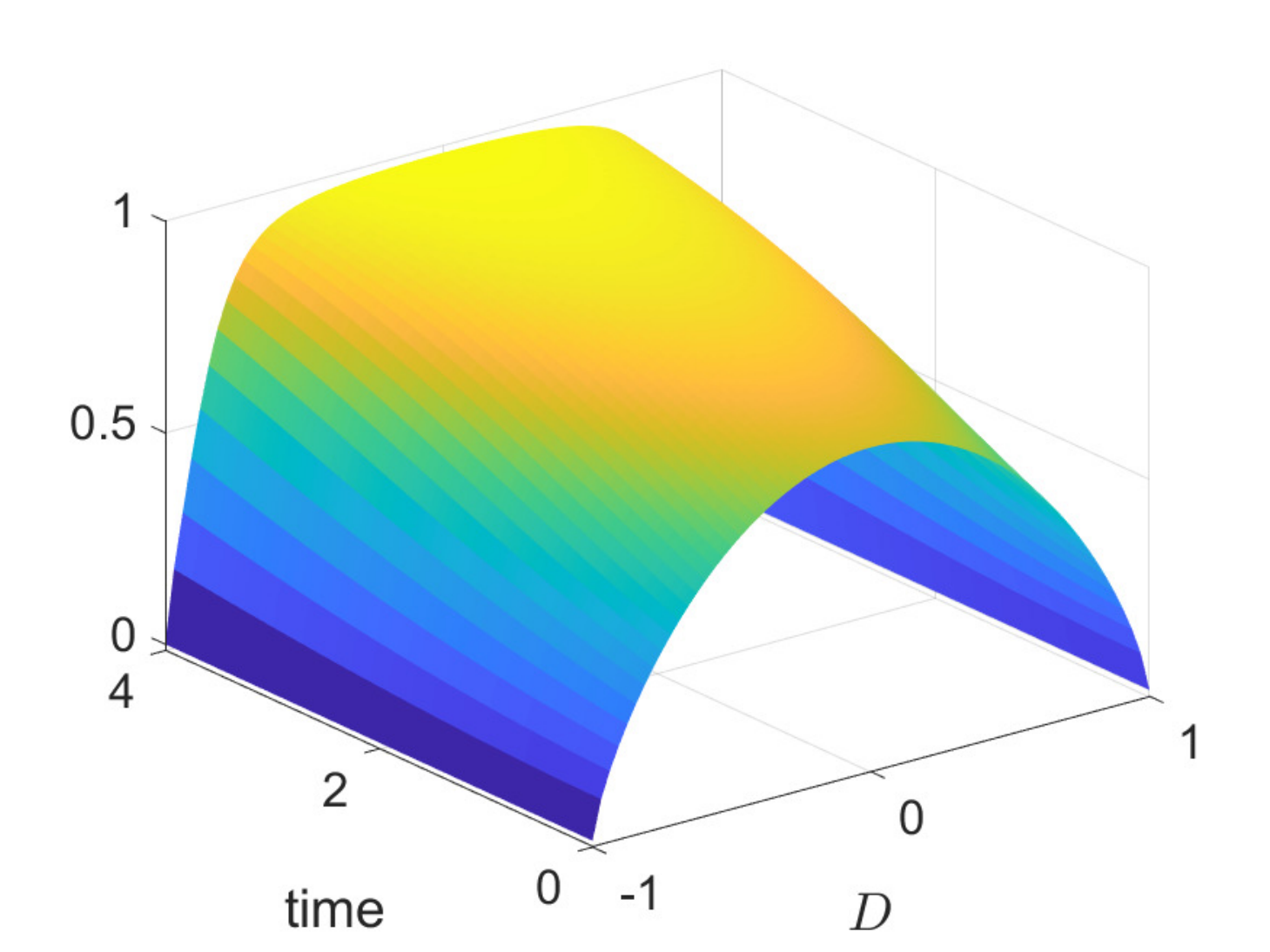}
	\includegraphics[scale = 0.39]{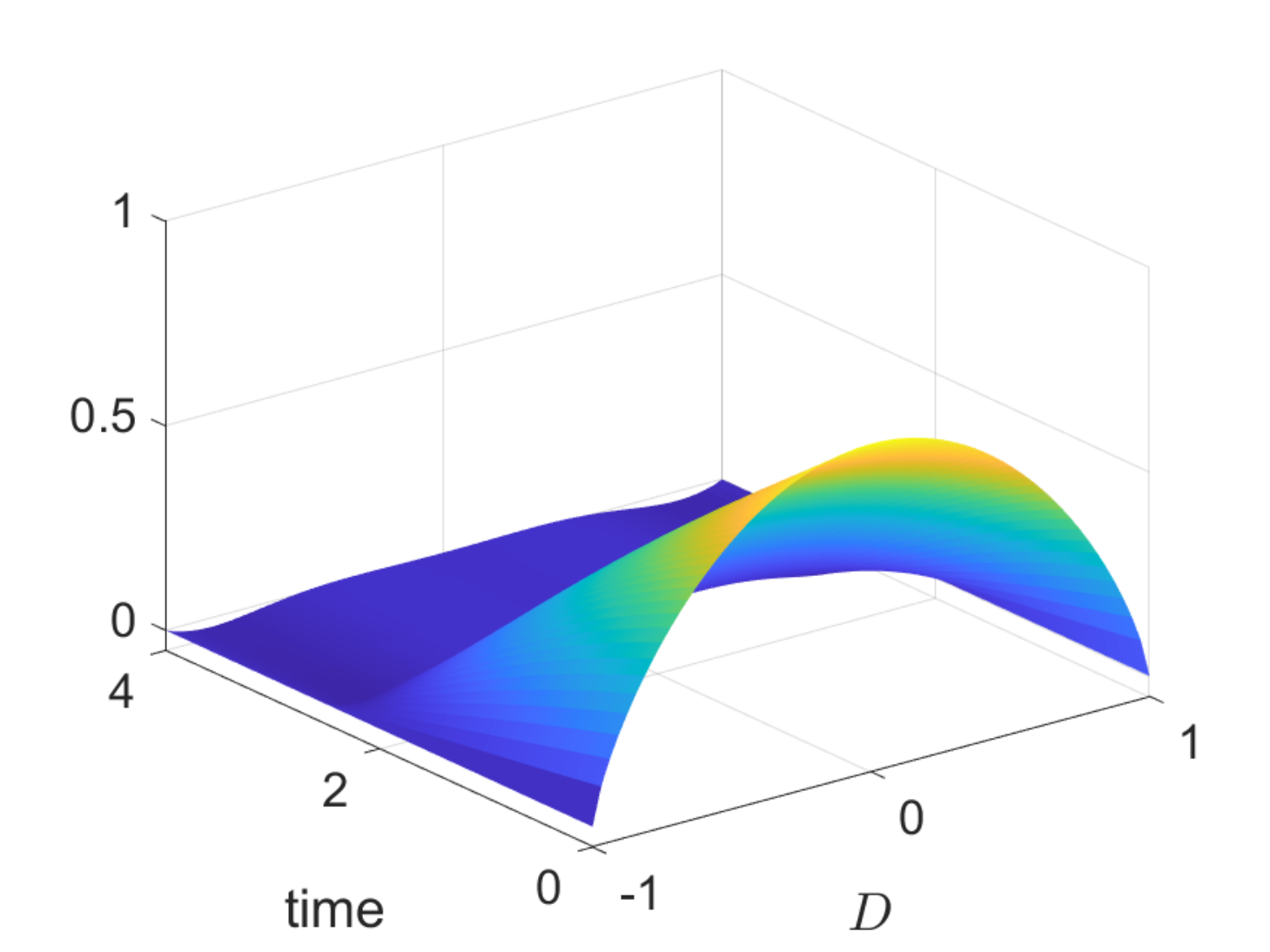}
	\caption{Test 3. Left: Uncontrolled solution. Right: Controlled solution.}
	\label{fig1_test3}
\end{figure}

\begin{figure}[htbp]
	\centering
	\includegraphics[scale = 0.39]{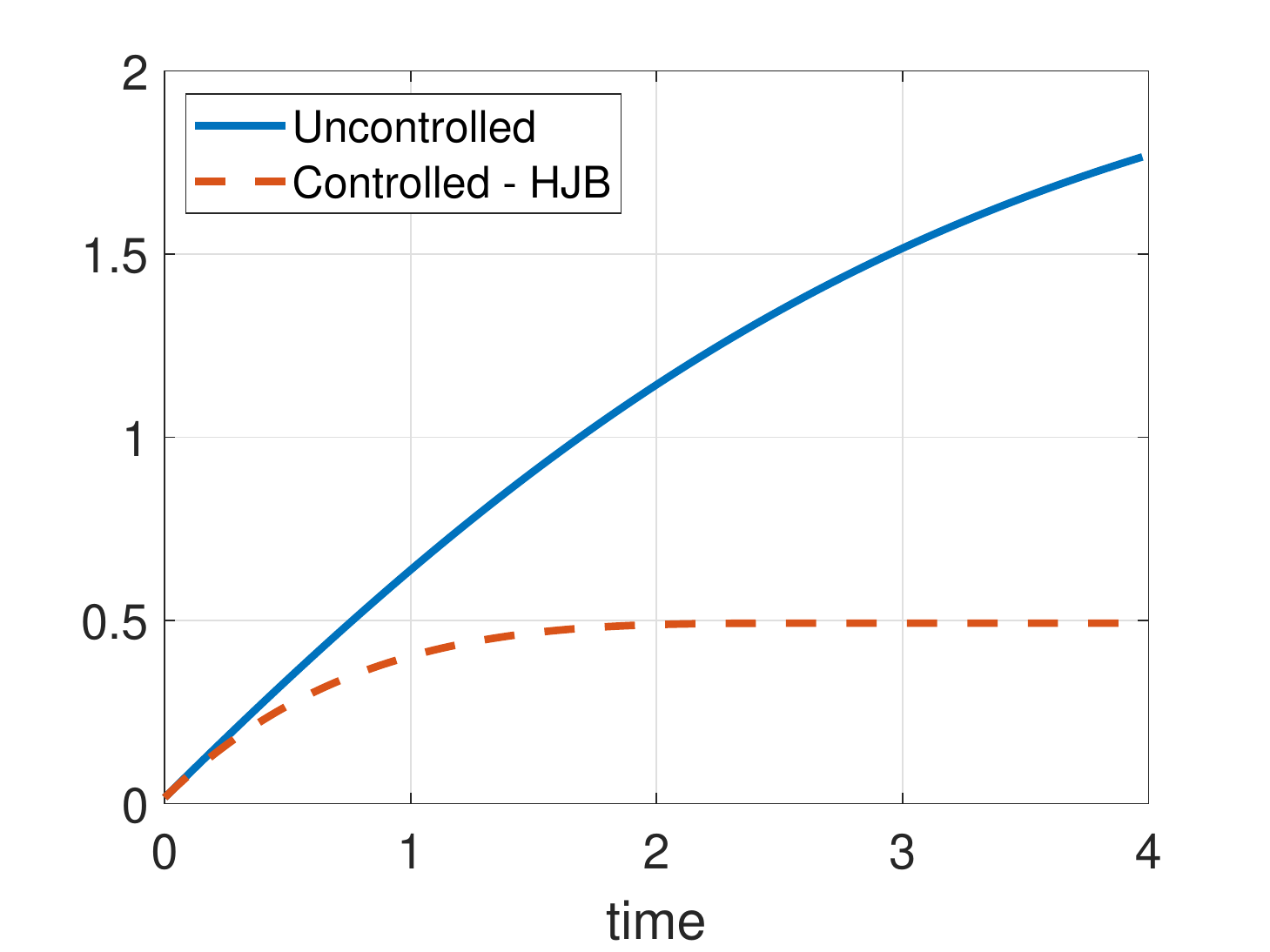}
	\includegraphics[scale = 0.39]{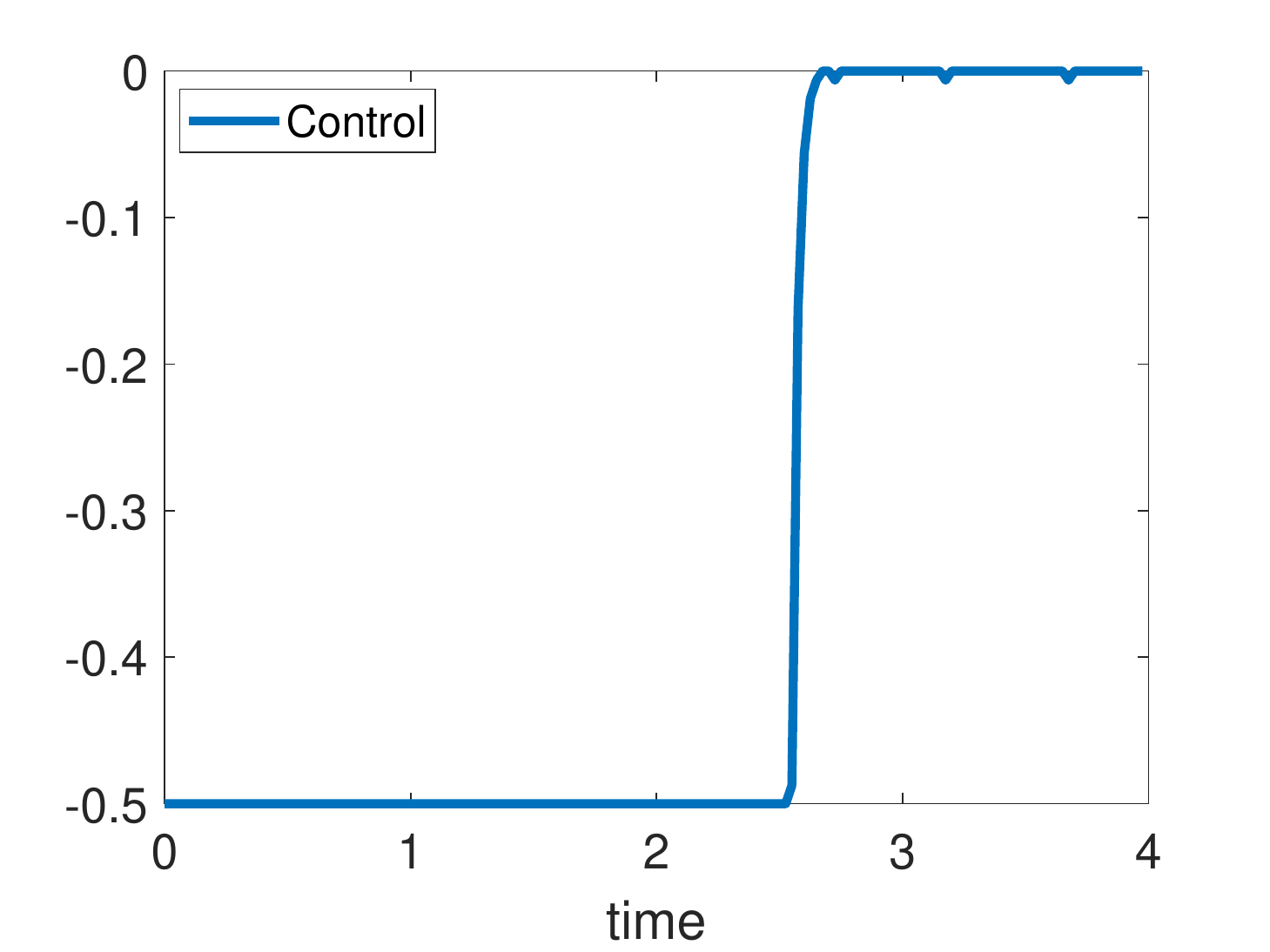}\\
	\caption{Test 3. Left: Cost functional. Right: Optimal control. }
	\label{fig2_test3}
\end{figure}

The space $\mathcal{P}=[0.08, 0.12]$ is discretized with step size $0.01$. The residual is minimized when the parameter is equal to $0.09$. We run the VI algorithm with $\Delta t = 0.01$ and $U$ discretized in $21$ nodes. We obtain the feedback control and solution considering $81$ discrete control points and $\Delta t = 0.025$. 
In Figure \ref{fig1_test3}, we show the uncontrolled solution on the left and the controlled solution on the right. We can observe how the controlled solution reaches the desired configuration which is the zero equilibrium.
In the left panel of Figure \ref{fig2_test3} we compare the evaluation of the cost functional for the controlled and uncontrolled solutions. As expected, the controlled solution has a lower cost functional than the the uncontrolled one. Finally, in the right panel of Figure \ref{fig2_test3} we show the obtained control. The control is active at the beginning then it is zero when the desired configuration is reached.

\begin{figure}[htbp]
	\centering 	
	\includegraphics[scale = 0.39]{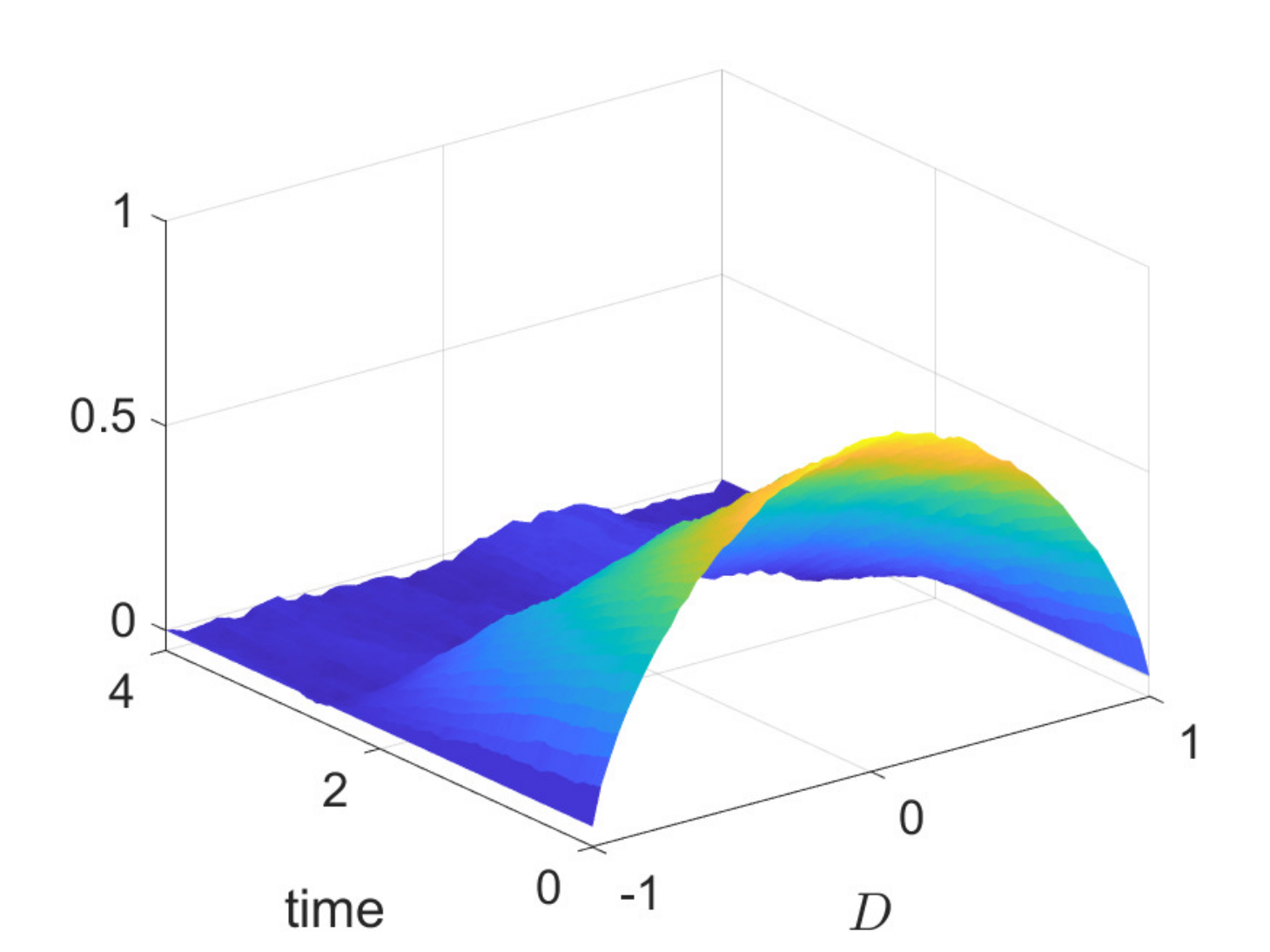}
	\includegraphics[scale = 0.39]{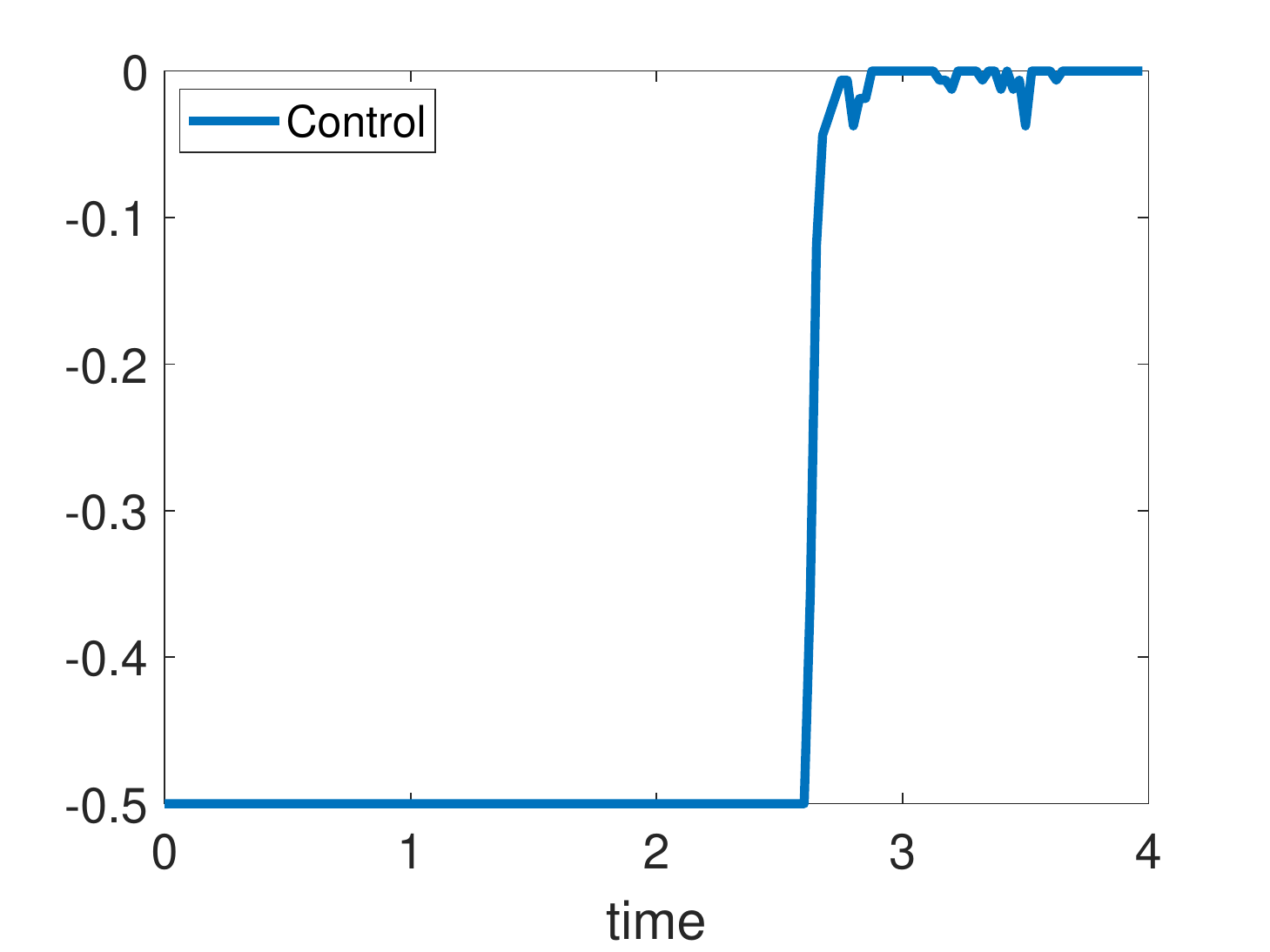}

	\caption{Test 3. Simulation with noise. Left: Controlled solution. Right: Optimal control.}
	\label{test3_noise}
\end{figure}

\paragraph{Simulation with noise}
Here we consider as perturbation a $63$ dimensional vector of independent, identically distributed Gaussian variables with zero mean and standard deviation $0.0025$.
At each instance of time we add to the trajectory a new independent perturbation term.
The feedback control and the trajectory are computed with $\Delta t = 0.025$ and $U$ discretized in $81$ nodes.
In the left panel of Figure \ref{test3_noise} we show the controlled solution under disturbances. We can see that, despite the noise, we are able to reach the desired configuration. Then, in the right panel of Figure \ref{test3_noise} we show the control found. The behaviour of the control is qualitatively similar to the picture in the right panel of Figure \ref{fig2_test3}, but the feedback is able to react and adjust, see e.g. the correction due to noise near t=3.5.


\section{Conclusions}\label{sec:conclusion}
In this paper we presented the first approach to the feedback control of nonlocal, fractional problems using the dynamic programming principle. 
A key feature of our approach is the combination of semi-Lagrangian schemes with Shepard approximation. The results discussed in the numerical section illustrate the accuracy of the method
and the ability to tackle nonlinearities and noise. As this is the first work in this direction, we would like to further investigate feedback control for different nonlocal models and applications.
Furthermore, a natural extension of this work includes considering different approaches such as Model Predictive Control \cite{GP17}; it would permit us to increase the resolution of the discretization of the state equation and improve the usability of feedback control for nonlocal problems.

\section*{Acknowledgements}
Sandia National Laboratories is a multi-mission laboratory managed and operated by National Technology and Engineering Solutions of Sandia, LLC., a wholly owned subsidiary of Honeywell International, Inc., for the U.S. Department of Energy’s National Nuclear Security Administration under contract DE-NA0003525. This paper, SAND2022-14035 O, describes objective technical results and analysis. Any subjective views or opinions that might be expressed in the paper do not necessarily represent the views of the U.S. Department of Energy or the United States Government.

The work of M. D'Elia and C. Glusa is supported by by the Sandia National Laboratories Laboratory Directed Research and Development (LDRD) program. M. D'Elia is also partially supported by the U.S. Department of Energy, Office of Advanced Scientific Computing Research under the Collaboratory on Mathematics and Physics-Informed Learning Machines for Multiscale and Multiphysics Problems (PhILMs) project.

The work of H. Oliveira was supported by the Coordenação de Aperfeiçoamento de Pessoal de Nível Superior - Brasil (CAPES) - Finance Code 001.\\

\appendix

\section{Derivation of example with analytic solution}
\label{sec:deriv-example}

Let \(q\) and \(\tilde{b}\) be the exact solution and corresponding right-hand side of the fractional Poisson problem
\begin{equation*}
	\left\{
	\begin{array}{rcll}
		 (-\Delta)^{\fracOrder}q(\xib) &=& \tilde{b}(\xib) &\xib\in \omg,  \\
		q(\xib)&=&0 &\xib\in \omg^{c}, \\
	\end{array}
	\right.
\end{equation*}
with \(\|q\|_{L^{2}(\omg)}=1\).
An example of such a pair when \(\omg=(-1,1)\) is given by
\begin{align*}
  \tilde{b}(\xib) &= 2^{2\fracOrder}\Gamma\left(1+\fracOrder\right) \frac{\Gamma(\fracOrder+1/2)}{\Gamma(1/2)} \sqrt{\frac{\Gamma(2\fracOrder+3/2)}{\Gamma(2\fracOrder+1)\Gamma(1/2)}}, \\
  q(\xib) &= \sqrt{\frac{\Gamma(2\fracOrder+3/2)}{\Gamma(2\fracOrder+1)\Gamma(1/2)}} \left(1-\xib^{2}\right)^{\fracOrder}_{+}.
\end{align*}
Now, let \(\phi,\kappa:(0,T)\rightarrow\mathbb{R}\) such that \(\phi(0)=1\), \(\kappa(T)=0\).
Let \(U:=[a,b]\) with \(a<0<b\) and let \(\gamma,\lambda>0\).
Consider the following functions that will be used in the construction of the cost functional and the state equation
\begin{align*}
	y_{d}(\xib,t) &:=\phi(t)q(\xib) - \gamma \kappa'(t)q(\xib) + \gamma \kappa(t) \tilde{b}(\xib) + \lambda\gamma\kappa(t)q(\xib), \\
	u_{d}(t) &:= \operatorname{proj}_{U}(\kappa(t)) ,\\
	b(\xib,t) &:= \phi'(t)q(\xib) + \phi(t)\tilde{b}(\xib) - u_{d}(t)q(\xib).
\end{align*}
We first aim to minimize the finite time horizon cost functional
\begin{align*}
	\mathcal{J}_{q}^{T}(y, u)
	&:= \frac{1}{2} \int_0^{T} (\|y(\cdot,\eta)-y_{d}(\cdot,\eta)\|_{L^{2}(\omg)}^{2}  + \gamma \|u(\cdot, \eta)\|_{L^{2}(\omg)}^2) e^{-\lambda \eta} d\eta \\
	&= \frac{1}{2} \|y-y_{d}\|_{L_{\nu}^{2}(0,T;\omg)}^{2} + \frac{\gamma}{2} \|u\|_{L_{\nu}^{2}(0,T;\omg)}^{2},
\end{align*}
where we have set \(\nu(t):=e^{-\lambda t}\),
subject to
\begin{equation}
	\label{heat_test1}
	\left\{
	\begin{array}{rll}
		\partial_t y(\xib,t) + (-\Delta)^{\fracOrder}y(\xib,t) &= b(\xib,t) + u(t)q(\xib) &(\xib,t) \in \omg \times (0,T){,}  \\
		y(\xib,t)&=0 &(\xib,t)\in \omg^{c} \times (0,T) , \\
		y(\xib,0)&=q(\xib) &\xib\in \omg .
	\end{array}
	\right.
\end{equation}

For fixed initial condition \(q\), we can write the first order optimality conditions.
Let \(j(u):=\mathcal{J}_{q}^{T}(\mathcal{S}u,u)\), where \(\mathcal{S}\) is the solution operator of the state equation \eqref{heat_test1} above. The operator $\mathcal{S}$ is also called the fractional control-to-state operator $$\mathcal{S}: L^{2}(0,T;\omg) \rightarrow \mathbb{V}$$
and $\mathcal{S}u=z(u)$, where $z(u)$ solves \eqref{heat_test1}.

Then, the optimal control must satisfy the variational inequality 
\begin{align*}
	\overline{u}=\operatorname{argmin} j(u) \Leftrightarrow ( j'(\overline{u}), u-\overline{u})\geq 0 \quad \forall u \in \mathcal{U}_{ad}.
\end{align*}

The inequality can be written in the equivalent form
\begin{align*}
	(\mathcal{S}^{*}(\mathcal{S}\overline{u}-y_{d}) + \gamma q\overline{u},qu-q\overline{u})_{L_{\nu}^{2}((0,T);\omg)}\geq0,
\end{align*}
where \(\mathcal{S}^{*}\) is the adjoint solution operator. Hence, \(\overline{p}:=\mathcal{S}^{*}(\mathcal{S}\overline{u}-y_{d})\) is the solution to
\begin{align*}
	\left\{
	\begin{array}{rll}
		-\partial_t \overline{p}(\xib,t) +\lambda \overline{p} + (-\Delta)^{\fracOrder}\overline{p}(\xib,t) &= \overline{y}(\xib,t)-y_{d}(\xib,t) &(\xib,t) \in \omg \times (0,T), \\
		\overline{p}(\xib,t)&=0 &(\xib,t)\in \omg^{c} \times (0,T), \\
		\overline{p}(\xib,T)&=0 &\text{for all }\xib\in \omg.
	\end{array}
	\right.
\end{align*}
with $\overline{z} = \mathcal{S}\overline{u}$.
Then, the following inequality
\begin{align*}
	0\leq(\overline{p} + \gamma q\overline{u},qu-q\overline{u})_{L_{\nu}^{2}(0,T;\omg)} = ((\overline{p},q)_{L^{2}(\omg)}q + \gamma q\overline{u},qu-q\overline{u})_{L_{\nu}^{2}(0,T;\omg)}
\end{align*}
implies that 	$\overline{u} = \operatorname{proj}_{U}\left(-\frac{1}{\gamma}(\overline{p},q)_{L^{2}(\omg)}\right)$.
If we set
\begin{align*}
	y^{*}(\xib,t):=\phi(t)q(\xib), \quad
	p^{*}(\xib,t):=-\gamma\kappa(t)q(\xib), \quad
	u^{*}(t):=u_{d}(t),
\end{align*}
we obtain
\begin{align*}
	y^{*}(\xib,0)
	&= q(\xib), \\
	\partial_t y^{*}(\xib,t) + (-\Delta)^{\fracOrder}y^{*}(\xib,t)
	&= \phi'(t)q(\xib) + \phi(t)\tilde{b}(\xib)
	= b(\xib,t)+u^{*}(\xib,t),
\end{align*}
where we have used that $q$ is the solution of the fractional Poisson equation. Hence, \(y^{*}\) is the state corresponding to the control \(u^{*}\) and the initial condition \(q\).
Moreover,
\begin{align*}
	p^{*}(\xib,T) &=0,\\
	-\partial_t p^{*}(\xib,t) + \lambda p^{*}(\xib,t) + (-\Delta)^{\fracOrder}p^{*}(\xib,t) &= \gamma\kappa'(t)q(\xib) -\lambda\gamma\kappa(t)q(\xib) -\gamma\kappa(t)\tilde{b}(\xib) \\
	&= y^{*}(\xib,t)-y_{d}(\xib,t).
\end{align*}
and hence \(p^{*}\) solves the adjoint equation with right-hand side \(y^{*}-y_{d}\).
Finally, $\operatorname{proj}_{U}\left(-\frac{1}{\gamma}(p^{*},q)_{L^{2}}\right) = \operatorname{proj}_{U}\left(\kappa(t)\right) = u^{*}(t).$
Therefore, for the initial condition \(q\), the optimal control is \(u^{*}\), and the optimal state is \(y^{*}\).

Let us now link the problem to the infinite time horizon framework. In order to do so, let us choose \(T_{0}>0\) and \(\kappa\) such that \(\kappa(t)=0\) for \(t\geq T_{0}\).
For \(T>T_{0}\), the previous construction gives the solution of the optimal control problem on \((0,T)\).
On the other hand, we have
\begin{align*}
	\mathcal{J}_{q}^{\infty}(y^{*},u^{*})
	&= \mathcal{J}_{q}^{T}(y^{*},u^{*}) + \frac{1}{2} \|y^{*}-y_{d}\|_{L_{\nu}^{2}(T,\infty;D)}^{2} + \frac{\gamma}{2} \|u^{*}\|_{L_{\nu}^{2}(T,\infty;D)}^{2} \\
	&= \mathcal{J}_{q}^{T}(y^{*},u^{*}) + \frac{\gamma^{2}}{2} \left[ \|\kappa'\|_{L_{\nu}^{2}(T,\infty)}^{2} + \|\kappa\|_{L_{\nu}^{2}(T,\infty)}^{2}\|\tilde{b}\|_{L^{2}(D)}^{2}\right. -\lambda (\kappa',\kappa)_{L_{\nu}^{2}(T,\infty)}\\
	&-\left.2(\kappa',\kappa)_{L_{\nu}^{2}(T,\infty)} (q,\tilde{b})_{L^{2}(D)} +2\frac{\lambda}{\gamma}\|\kappa\|_{L_{\nu}^{2}(T,\infty)}^{2}(q,\tilde{b})_{L^{2}(D)} + \lambda^2  \|\kappa\|_{L_{\nu}^{2}(T,\infty)}\right] \\ & \qquad + \frac{\gamma}{2} \|\operatorname{proj}_{U}(\kappa)\|_{L_{\nu}^{2}(T,\infty)}^{2} \\
	&= \mathcal{J}_{q}^{T}(y^{*},u^{*}),
\end{align*}
where we have used that \(\kappa\) and \(\kappa'\) are zero on \((T,\infty)\). Therefore,
\begin{align*}
	\min \mathcal{J}_{q}^{\infty}(y,u) \leq \min \mathcal{J}_{q}^{T}(y,u),
\end{align*}
but the inverse inequality also holds, since \(\mathcal{J}_{q}^{T}(y,u)\le\mathcal{J}_{q}^{\infty}(y,u)\).
Hence, the pair \((y^{*},u^{*})\) is also optimal for the infinite time horizon case.


\section*{Declarations}
 
\paragraph{Ethical Approval}
not applicable
 
\paragraph{Competing interests}
No conflict of interest.

\paragraph{Authors' contributions}
A.A: idea of applying to nonlocal, writing, numerical tests. M.D. : writing, funding. C.G.: fem simulations, writing. H.O.: numerical tests, figures.
 
\paragraph{Funding}
MD anc CG were financed by LDRD. HO by CAPES
 
\paragraph{Availability of data and materials} 
not applicable

\bibliographystyle{plain}
\bibliography{snl,hjb}
\end{document}